\documentclass[11pt]{amsart}
\usepackage{amsmath,amsfonts,amssymb,amsxtra,latexsym,amscd,enumerate,amsthm}

\usepackage{latexsym}
\usepackage{epsfig}

\def\z{\zeta}
\def\t{\theta}

\def\a{\alpha}
\def\q{\frac{1}{2}}

\def\e{\varepsilon}

\def\D{\Delta}

\def\Ra{\Rightarrow}
\def\les{\lesssim}

\def\R{\mathbb{R}}
\def\C{\mathbb{C}}
\def\Z{\mathbb{Z}}
\def\N{\mathbb{N}}
\def\S{\mathbb{S}}

\def\beq{\begin{equation}}
\def\eeq{\end{equation}}
\def\eq{\Leftrightarrow}

\def\beq{\begin{equation}}
\def\eeq{\end{equation}}
\def\ra{\rightarrow}
\def\Leq{\Longleftrightarrow}

\newtheorem{t1}{Theorem}
\newtheorem{l1}{Lemma}
\newtheorem{p1}{Proposition}

\newtheorem*{r2}{Remark}

\begin{document}
\title[On Schr\"odinger Maps]{On Schr\"odinger Maps}

\author{Ioan Bejenaru}
\address{Department of Mathematics, UCLA, Los Angeles CA 90095-1555}
\email{bejenaru@math.ucla.edu}

\begin{abstract} 

We study the local well-posedness theory for the Schr\"odinger Maps equation. We work in $n+1$ dimensions, for $n \geq 2$, and prove a local well-posedness for small initial data in $H^{\frac{n}{2}+\e}$.

\end{abstract}

\maketitle

\section{Introduction}

If $(M,g)$ is a Riemannian manifold, then the harmonic maps are smooth maps $\phi:\R^{n} \rightarrow M$ which minimize the Lagrangian:

$$L_g^h=\q \int_{\R^n} |\nabla \phi|_g^2 dx$$

The corresponding Euler-Lagrange equations give rise to the so called harmonic map equation which is one of the most studied equations in the modern geometric analysis. 

If $\R^{n+1}$ is replaced by the Minkowski space, then the Lagrangian is adjusted in a natural way:

$$L_g^w=\q \int_{\R^n \times \R} -|\partial_t \phi|_g^2 +|\nabla_x \phi|_g^2 dx$$

\noindent
and the corresponding Euler-Lagrange equations are referred to as the wave maps equations.  The wave maps are maps $\phi : \R^{n+1} \rightarrow M$ which satisfy the wave map equations. We refer the reader to the expository work of Tataru, see \cite{ta2}, for a complete introduction in the field of the wave maps. From the same work we import the form of the wave map equations in local coordinates:

\beq \label{w}
\Box \phi^i=\Gamma^i_{jk}(\phi) \partial_{\a} \phi^{j} \partial_{\a} \phi^k
\eeq

\noindent
where $\Gamma^i_{jk}$ are the Christofell symbols and $\Box=\partial_{tt}-\Delta_{x}$. The first major breakthrough in treating the equation in the form (\ref{w}) was the observation due to Klainerman, see \cite{kl}, that the nonlinearity is not generic. Klainerman noticed that the derivative part of the nonlinearity, which we can simplified to $(\nabla{\phi})^2$, exhibits a special structure which is called the null condition. Since then there has been a lot of development in the subject, see Tataru \cite{ta1} and Tao \cite{tao1} for most recent results on the wave maps. 

Their Schr\"odinger counterpart, the Schr\"odinger maps, has been studied only recently. If we assume that for every positive or negative result for wave maps there should be an equivalent result for the the Schr\"odinger maps, then one could say that the theory for Schr\"odinger maps is well behind the one for the wave maps. We will make this more clear once we setup in a rigorous way the mathematical formulation of the problem. 

There are several ways to introduce the Schr\"odinger maps and they are all equivalent for smooth solutions. An elaborate introduction can be found in the work of Nahmod, Stefanov and Uhlenbeck, see \cite{n1}. 

We prefer to introduce them the same way we did with the harmonic and wave maps. Here and throughout the rest of the paper we choose $M=\S^2$ and we identify $\S^2 \setminus \N$ ($\N$ is the north pole) with the Riemannian surface $(\C, g dz d \bar{z})$ by using the stereographic projection:

$$z \in \C \mapsto (\frac{2 \mbox{Re} z}{1+|z|^2}, \frac{2 \mbox{Im} z}{1+|z|^2}, \frac{1-|z|^2}{1+|z|^2}) \in \S^2$$  

$g$ is given by $g(z,\bar{z})=(1+|z|^2)^{-2}$. By working with $H^{\frac{n}{2}+\e} \subset L^{\infty}$, we keep the solutions local and avoid the problematic issue that this representation has close to the north pole.  For each $t \in \R$, the energy of the map $z: \R^n \times \R \rightarrow \mathbb{C}$ is defined by:

\beq \label{e1}
L^s_g(z(t))=\q \int_{\R^n} \frac{|\nabla z|^2}{(1+|z|^2)^2} dx
\eeq

The Euler-Lagrange equations of this energy functional are:

\beq \label{e2}
\sum_{j=1}^{n} (\frac{\partial}{\partial x_j}- \frac{2 \bar{z}}{1+|z|^2} \frac{\partial z}{\partial x_j}) \frac{\partial z}{\partial x_j}=0
\eeq

Since the pullback covariant derivative by the map $z$ from $\R^n \times \R$ to $(\C, g dz d \bar{z})$ is given by:

\beq \label{e3}
\nabla_j=\frac{\partial}{\partial x_j}- \frac{2 \bar{z}}{1+|z|^2} \frac{\partial z}{\partial x_j}
\eeq

\noindent
the expression on the left in (\ref{e2}) is the most natural Laplacian on $(\C, g dz d \bar{z})$. The Schr\"odinger maps equation is the following evolution equation:

\beq \label{e4}
i\frac{\partial z}{\partial t}= \sum_{j=1}^{n} (\frac{\partial}{\partial x_j}- \frac{2 \bar{z}}{1+|z|^2} \frac{\partial z}{\partial x_j}) \frac{\partial z}{\partial x_j}
\eeq

We can rewrite the equation as:

\beq \label{e5}
i z_{t} - \Delta z=  \frac{2 \bar{z}}{1+|z|^2} (\nabla z)^2
\eeq

\noindent
and we obtain a representation similar to the one in (\ref{w}). The equation (\ref{e5}) is a nonlinear Schr\"odinger equation with derivatives (D-NLS). Among people working on Schr\"odinger maps there is almost a general consensus that it is not lucrative to exploit the equation in the form (\ref{e5}) and there are several reasons for not favoring a D-NLS. One of them is that the standard energy methods cannot be applied in order to derive a local well-posedness result. The other one is that even more sophisticated methods (see Kenig-Ponce-Vega \cite{k2}, Chihara \cite{ch1}) require a high level of regularity on the initial data and some decay structure. A more recent result of the author, see \cite{be}, succeeds to lower the level of regularity all the way down to the critical exponent for generic D-NLS in two dimensions. But the initial data has to have a bit of spherical symmetry and some decay at infinity and these are not natural conditions to impose in the case of the Schr\"odinger maps equation. 

Therefore the general approach for dealing with the Schr\"odinger maps equation was to introduce a gauge transform which creates a moving frame depending on the solution that puts the pullback of the covariant derivative $\nabla_j$ as near to $\frac{\partial}{\partial x_j}$ as possible. This gives rise to a new systems of equations which is called the modified Schr\"odinger maps (MSM). The system contains the equations in the new coordinates plus an additional system describing the evolution of the moving frame. We do not follow this approach and this is why we skip the technical details of the construction; for reference the reader can consult \cite{n1}.  

In order to understand better what has been achieved so far in the subject, we describe now what are the expectations from the Schr\"odinger maps equation. The scaling exponent of the problem is $s_c=\frac{n}{2}$ and the solutions of (\ref{e5}) preserve the energy $L^s_g$ in (\ref{e1}). Following the heuristics in \cite{ta2} one expects the following kind of development of the theory:

\vspace{.05in}

\begin{itemize}

\item if $s>\frac{n}{2}$, establish local well-posedness for initial data $u_0 \in H^s$

\vspace{.05in}

\item if $s=\frac{n}{2}$, establish global well-posedness for small initial data in $u_0 \in \dot{H}^s$

\vspace{.05in}

\item if $s=\frac{n}{2}$, analyze whether the large energy solutions blow-up in finite time.

\end{itemize}

The best results known so far are as follows. In the two dimensional case, Kato \cite{ka} and, independently, Kenig and Nahmod \cite{k1} showed existence of (MSM) with initial data in $H^{\q+\e}$. Later, Kato and Koch \cite{he} showed uniqueness for (MSM) with initial data in $H^{\frac{3}{4}+\e}$. In order to read these results correctly we should note that the gauge transformation giving (MSM) involves a derivative of the solution to the Schr\"odinger maps equation. Therefore a result at the level of $H^s$ for the (MSM) corresponds to a result at the level of $H^{s+1}$ for the Schr\"odinger maps. In other words, for the original Schr\"odinger maps equation in two dimension we know existence of solutions for initial data in $H^{1+\q+\e}$ and uniqueness for initial data in $H^{1+\frac{3}{4}+\e}$.

We learned very recently, by personal communication, that A. Ionescu and C. Kenig have just obtained a new result of the same nature we state in this paper. We will update our paper, with the necessary commentaries, once a copy of their paper is available. Our work and Ionescu and Kenig's work have been carried out completely independently.

Our goal in this paper is to provide a local well-posedness result for the Schr\"odinger maps equation with initial data in $H^{\frac{n}{2}+\e}$. We do not follow the standard approach, i.e. solve the problem via the (MSM), but work directly with the equation in the form (\ref{e5}). We are aware of the fact that previous theories on the generic D-NLS required additional structure(s) on the initial data especially for the low regularity case. 

Our reason to treat the Schr\"odinger maps equation as a D-NLS comes from the same kind of observation that Klainerman made for the wave maps: we do have a non-generic D-NLS. The derivative part of the nonlinearity, $(\nabla{z})^2$, has a special structure in it which kills the worst interactions in D-NLS. We try to sketch an argument for this claim. 

If $\phi$ is a homogeneous solution of the Schr\"odinger equation equation then its Fourier transform $\hat{\phi}$ is supported on the paraboloid $P=\{(\xi,\tau) \in \R^n \times \R: \tau=\xi^2 \}$. From this observation, we build our spaces such that the Fourier transform of their elements have a certain rate of decay away from $P$. When trying to iterate the derivative part of the nonlinearity, it turns out that the worst interactions are the ones coming from $P$ and ending on $P$. We formalize this bellow. If $\phi^1$ and $\phi^2$ are two functions in our space, the derivative part of their interaction is

$$\widehat{\nabla \phi_1 \nabla \phi_2}(\xi,\tau)= \int  \widehat{\nabla{\phi_1}}(\xi^1,\tau^1) \widehat{\nabla{\phi_2}}(\xi^2,\tau^2)  d \xi^1 d\tau^1=$$

$$- \int  \xi^1 \cdot \xi^2 \hat{\phi}_1(\xi^1,\tau^1) \hat{\phi}_2(\xi^2,\tau^2)  d \xi^1 d\tau^1$$

\noindent
with the convention that $\xi=\xi^1+\xi^2$ and $\tau=\tau^1+\tau^2$. The worst interactions are those of type $(\xi^1, (\xi^1)^2)$ and $(\xi^2, (\xi^2)^2)$ with the property that $(\xi^1+\xi^2, (\xi^1)^2+(\xi^2)^2) \in P$. This is equivalent to saying that $\xi^1 \cdot \xi^2=0$. Therefore the worst interactions are absent in the above computation and this is why we say that the derivative part of the nonlinearity has a killing structure similar to the one for wave maps. This does not make the problem trivial. But combined with previous experience of the author with D-NLS, see \cite{be}, it led us to a different perspective on the subject and a new result.

We describe bellow our approach in solving the Schr\"odinger maps equation in its local coordinate form (\ref{e5}). We look at the more general semilinear equation:

\begin{equation} \label{E}
\begin{cases}
\begin{aligned} 
iu_{t}-\Delta u &= Q(u,\bar{u}) (\nabla u)^2, \   t\in \mathbb{R}, x\in \mathbb{R}^{n} \\
u(x,0) &=u_{0}(x)
\end{aligned}
\end{cases}
\end{equation}

\noindent
where $u: \mathbb{R}^{n} \times \mathbb{R} \rightarrow \mathbb{C}$ and $Q: \C \times \C \rightarrow \C$ is analytic. We denote by $N(u)=Q(u,\bar{u}) (\nabla u)^2$ the nonlinearity of our problem. Our nonlinearity in (\ref{e5}) can be written this way as long as we keep $||u||_{L^{\infty}}$ small and we will make sure this is the case.

We introduce the corresponding inhomogeneous Schr\"odinger equation:

\begin{equation} \label{I}
\begin{cases}
\begin{aligned} 
iu_{t}-\Delta u &= f, \   t\in \mathbb{R}, x\in \mathbb{R}^{n} \\
u(x,0) &=g(x)
\end{aligned}
\end{cases}
\end{equation}

We seek for spaces of functions $Z^s$ and $W^s$ with the following properties:

\begin{itemize}

\item{(\it linear property)} \rm The solution $u$ to (\ref{I}) satisfies

\beq \label{le}
||\chi(t) u||_{Z^{s}} \les ||g||_{H^s}+||f||_{W^s}
\eeq

\item{ \it (nonlinear estimate)} \rm The truncated $\chi(t)N$ has the mapping property

$$\chi(t) N: Z^s \rightarrow W^s \ \ \mbox{is Lipschitz continuous}$$

\end{itemize}

The most difficult part in these kind of problems is to construct the spaces the right way and to be able to prove the nonlinear estimates. Proving the nonlinear estimates involves solving two problems. We know from the theory of wave maps that $Z^s$ should be an algebra since this way $Q(u)$ (we deliberately ignore that $Q$ involves $\bar{u}$ too) maps $Z^s$ to $Z^s$ when $||u||_{Z^s}$ is small. Our problem involves functions with complex values and as we will see in the next section, $Z^s$ is not closed under conjugation. For a space $X$ we define its conjugate $\bar{X}$ by:

\beq \label{C}
u \in \bar{X} \eq \bar{u} \in X
\eeq

Therefore we will have to show that 

\beq \label{A}
Z^s+\bar{Z}^s \ \ \mbox{is an algebra}
\eeq

\noindent 
in order to conclude that if $u \in Z^s$ with small norm, then $Q(u,\bar{u}) \in Z^{s}+\bar{Z}^s$. 
The next step is to show that $Z^s+\bar{Z}^s$ leaves $W^s$ invariant under multiplication:

\beq \label{M}
(Z^s+\bar{Z}^s) W^s \subset W^s
\eeq

Finally we need the bilinear estimate. If we denote by $B(u,v)=\nabla{u} \cdot \nabla{v}$ we have to show that:

\beq \label{be}
||B(u,v)||_{W^{s}} \leq  C_s ||u||_{Z^s} ||v||_{Z^s}
\eeq

Proving the nonlinear estimate amounts to proving (\ref{A})-(\ref{be}). One would notice that we do not really need time truncation on the bilinear term, as written in statement of the nonlinear estimate.

Once we have the linear property and the nonlinear estimate, then a standard fixed point argument solves the problem:

\begin{equation} 
\begin{cases}
\begin{aligned} 
iu_{t}-\Delta u &= \chi(t) N(u), \\
u(x,0) &=u_{0}(x)
\end{aligned}
\end{cases}
\end{equation}

If $\chi(t)$ is $1$ on $(-T,T)$ then we obtain that $u$ solves (\ref{E}) on the interval $(-T,T)$ an we can claim

\begin{t1} \label{T}

For any $n \geq 2$ and $s > \frac{n}{2}$ and $T > 0$, there exists $\delta > 0$ such that for every $u_{0} \in H^{s}$ with $||u_{0}||_{H^{s}}  < \delta$ , the problem (\ref{E}) has a unique solution $u \in C([0,T]:H^{s}) \cap Z^{s}$ with Lipschtiz dependence on the initial data.

\end{t1}

Theorem \ref{T} is a small data type result and there are two reasons for this limitation. One is that we use the analyticity of $Q$ and this imposes a smallness on the $||u||_{L^{\infty}}$ which is guaranteed by the smallness of $||u||_{Z^s}$. The second reason is that, even if $Q$ were analytic everywhere, we cannot just rescale the problem and get the large data result at the expense of shrinking the life-span of the solution. In other words, dealing with large data is a different problem and we did not find a trivial way to derive it from our result. This is probably the only downside of our result versus the set of all the other results in the field.

On the bright side, we should remark that our Theorem comes with the whole package of properties one would expect from a local well-posedness result. We mean existence, uniqueness and Lipschitz dependence on the initial data. In the previous results these were treated as different problems.

The paper is organized as follows. In the next section we introduce the notations and definitions we use along the paper. Section \ref{les} is devoted to the linear estimates while sections \ref{bes} and \ref{bes+} are devoted to the bilinear estimates. We end the paper with section \ref{alg} in which we establish the algebra type properties, see (\ref{A})-(\ref{M}). 

Our spaces and proofs are adapted for the case $T=1$. By rescaling one can easily deduce the result in Theorem \ref{T} for all $T >0$. 

\vspace{.1in}

The author thanks Daniel Tataru and Terry Tao for useful discussions and encouragement with the project, and to Herbert Koch and Andrea Nahmod for help on various aspects of the problem. 

\section{Notations and Definitions} \label{ds}

The notations in this paper are very similar (maybe with just one notable difference) to the ones the author has used in \cite{be}.

For each $u$ we denote by $\mathcal{F}u=\hat{u}$ the Fourier transform of $u$. This is always taken with respect to all the variables, unless otherwise specified.

Throughout the paper $A \lesssim B$ means $A \leq CB$ for some constant $C$ which is independent of any possible variable in our problem. We say $A \approx B$ if $A \leq CB \leq C^{2}A$ for the same constant $C$. We say that we localize at frequency $2^{i}$ to mean that  in the support of the localized function $\langle (\xi,\tau) \rangle \in [2^{i-1},2^{i+1}]$.

In the Schr\"odinger equation time and space scale in a different way, and this suggests to define the norm for $(\xi, \tau)$ by $|(\xi, \tau)|=(|\tau|+ \xi^{2})^{\q}$. We introduce the Bourgain space $X^{s,b}$:

$$X^{s,b} =  \{f\in S'; \langle (\xi,\tau) \rangle^{s} \langle \tau-\xi^{2} \rangle^{b} \hat{f} \in L^{2} \}$$

Here and thereafter $\langle x \rangle=(1+|x|^{2})^{\q}$ where $|x|$ is the norm of $x$. We employ frequency localized versions of $X^{s,\q}$ whose constructions are described bellow. 

Consider $s_{0} : [0,\infty) \rightarrow \R$ be a nonnegative smooth function such that $s_{0}(x)=1$ on $[0,1]$ and $s_{0}(x)=0$ if $x \geq 2$. Then for each $i \geq 1$ we define $s_{i}: [0,\infty) \rightarrow \R$ by $s_{i}(x)= s_{0}(2^{-i}x)- s_{0}(2^{-i+1}x)$. We define the operators $S_{i}$, to localize at frequency $2^{i}$, by:

$$\mathcal{F}(S_{i}f)= \hat{f}_{i} = s_{i}(|(\xi, \tau)|) \cdot \hat{f}(\xi, \tau)$$

We denote by $A_i$ the support of $s_i$ (more exactly of $s_i \circ |\cdot|$) in $\R^{n+1}$. The second type of localization we need is with respect to the modulation $\tau-\xi^2$. If $(\xi,\tau) \in A_{i}$, then $|\tau-\xi^{2}| \leq |\tau|+\xi^{2} \leq 2^{2i+2}$. For $d \in I_{i}=\{2^0, 2^1, .., 2^{2i+2}\}$ we define $s_{i,d}(\xi,\tau)=s_{i}(|(\xi, \tau)|) \cdot s_{\ln_{2}{d}}(|\tau-\xi^{2}|)$. The support of $s_{i,d}$ is denoted by $A_{i,d}$ and it can be described by:

$$A_{i,d}=\{(\xi,\tau): \langle (\xi,\tau) \rangle \approx 2^{i} \ \mbox{and} \  \langle \tau-\xi^2 \rangle \approx d \}$$

If $(\xi,\tau) \in A_{i,d}$ then $d \approx |\tau-\xi^{2}| \approx |(\tau,\xi)| d((\xi,\tau),P) \approx 2^{i} d((\xi,\tau),P)$ (away from zero) which implies $d((\xi,\tau),P) \approx 2^{-i}d$. It is important to have this estimate, since this gives us the size of $A_{i,d}$ in the normal direction to $P$ to be $\approx 2^{-i}d$. We have the property:

$$ \sum_{d \in I_{i}} \varphi_{i,d} (\xi,\tau) = \varphi_{i}((|\xi,\tau)|)  , \ \forall (\xi,\tau) \in \R^{n} \times \R $$

We define the operators $S_{i,d}$ by $S_{i,d}f=f_{i,d}=\check{s}_{i,d} * S_{i}f$ and: 

$$S_{i, \leq d}f = f_{i, \leq d}=\sum_{d' \in I_{i}: d' \leq d} f_{i, d'} \ \  \mbox{and} \ \ S_{i, \geq d}f=f_{i, \geq d}= f_{i}-f_{i, \leq d}$$

As before we can define $A_{i, \leq d}$ and $A_{i, \geq d}$ to be the support of the corresponding operators. The part of $\hat{f}$ which is at distance less than $1$ from $P$ plays an important role and this is why we define the global operators:

$$S_{P_{\leq 1}}f=f_{P_{\leq 1}}=\sum_{i=0}^{\infty} f_{i, \leq 2^i}  \ \mbox{and} \ S_{P_{\geq 1}}f=f_{P_{\geq 1}}=\sum_{i=0}^{\infty} f_{i, \geq 2^i}$$

For functions whose Fourier transform is supported in $A_{i}$  and for any $1 \leq p \leq \infty$ we define: 

$$||f||^{p}_{X_{i}^{0,\q,p}} = \sum_{d \in I_{i}} ||f_{i,d}||_{X^{0,\q}}^{p}$$

\noindent
with the usual convention for $p=\infty$. Then we define $X^{s,\q,p}$ by the norm:

$$||f||^{2}_{X^{s,\q,p}}=\sum_{i} 2^{2is} ||f_{i}||_{X^{0,\q,p}_{i}}$$

For technical purposes we need localized versions of this spaces, like $ X^{s,\q}_{i,d} = \{f \in X^{s,\q}: \hat{f} \ \mbox{supported in} \ A_{i,d} \}$ and, similarly, $ X^{s,\q,p}_{i,\leq d}$ and $ X^{s,\q,p}_{i,\geq d}$.

$X^{s,\q,1}$ is our first candidate for the space of solutions. Our computations indicate that it is the right space to measure only $f_{P_{\geq 1}}$ (theres is though a certain degree at freedom which will be made clear later); $f_{P_{\leq 1}}$ will be measured in $X^{s,\q,\infty}$ plus an additional structure whose construction is described bellow. 

In the plane $\tau=0$ we define the lattice $\Xi = \Z^n$. For each $\xi \in \Xi$ we build a non-negative function $\phi_{\xi}$ to be a smooth approximation of the characteristic function of the cube of size $\q$ in $\R^{2}$ centered at $\xi$ and satisfying the natural partition property:

$$\sum_{\xi \in \Xi} \phi_{\xi}=1$$

We can easily impose uniforms bounds on the derivatives of the system $(\phi_{\xi})_{\xi \in \Xi}$. For each $\xi \in \Xi$ with $|\xi| \approx 2^{i}$ we define:

$$f_{\xi} = \check{\phi}_{\xi} * f \ \ \ \ \mbox{and} \ \ \ f_{\xi, \leq d} = \check{\phi}_{\xi} * f_{i,\leq d}$$

The convolution above is performed with respect to the $x$ variable, i.e. it does not involve the $t$ variable.  

 We denote by $(Q^{m})_{m \in \Z^{n}}$ the standard partition of $\R^{n}$ in cubes of size $1$; i.e. $Q^{m}$ is centered at $m=(m_1,...,m_n) \in \Z^{n}$, has its sides parallel to the standard coordinate axis and has size $1$. For each $\xi \in \Xi$, $m \in \Z^{n}$ and $l \in \Z$ we define the tubes:

$$T_{\xi}^{m,l}=\cup_{t \in [l,l+1]} (Q^{m}-2t\xi) \times \{ t \}=$$

$$\{(x_1-2t\xi_{1},..., x_n-2t\xi_{n},t): (x_1,..,x_n) \in Q^{m} \ \mbox{and} \ t \in [l,l+1] \}$$

Then, for each $\xi \in \Xi$, we define the space $Y_{\xi}$ by the following norm:

$$||f||^{2}_{Y_{\xi}}= \sum_{(m,l) \in \Z^{n+1}} ||f||^{2}_{L^{\infty}_{t}L^{2}_{x}(T_{\xi}^{m,l})}$$

We have $f=\sum_{\xi \in \Xi} f_{\xi}$ and then we define the space $Y^{s}$ by the norm:

$$||f||^{2}_{Y^{s}}= \sum_{\xi \in \Xi} \langle \xi \rangle^{2s} ||f_{\xi}||^{2}_{Y_{\xi}}$$

We introduce the localized versions $Y_{i}= \{ f \in Y^{0}; \hat{f} \ \mbox{supported in} \ A_{i} \}$ and $Y_{i, \leq d}= \{ f \in Y^{0}; \hat{f} \ \mbox{supported in} \ A_{i, \leq d} \}$.

To bring everything together, define $Z^{s}$ to be 

$$Z^{s}=\{ f \in S': ||f_{P_{\geq 1}}||_{X^{s,\q,1}}+||f_{P_{\leq 1}}||_{Y^{s}}+||f_{P_{\leq 1}}||_{X^{s,\q,\infty}} < \infty \}$$

\noindent
with the obvious norm. 

We are done with the definition of $Z^{s}$, the space for our solutions, and continue with defining $W^s$. We can easily define $X^{s,-\q,p}$ by simply replacing $\q$ with $-\q$ in the definition of $X^{s,\q,p}$. Then we define $\mathcal{Y}_{\xi}$  and $\mathcal{Y}^{s}$ by:

$$||f||^{2}_{\mathcal{Y}_{\xi}}= \sum_{(m,l) \in Z^{n}} || f||^{2}_{L^{1}_{t}L^{2}_{x}(T_{\xi}^{m,l})}  \ \ \  \mbox{and} \ \ \  ||f||^{2}_{\mathcal{Y}^{s}}= \sum_{\xi \in \Xi} \langle \xi \rangle ^{2s} ||f_{\xi}||^{2}_{\mathcal{Y}_{\xi}}$$

Notice that $(\mathcal{Y}_{\xi})^{*}=Y_{\xi}$ since we will use this later for duality purposes.

 We introduce $W^{s}$ defined by the norm:

$$||f||_{\mathcal{W}^{s}}= \inf \{||f_{1}||_{\mathcal{Y}^{s}} + ||f_{2}||_{X^{s,-\q,1}}; f=f_{1}+f_{2} \}$$

We measure the right hand side of ($\ref{E}$) in: 

$$W^{s} =\{ f \in S': ||f_{P_{\leq 1}}||^{2}_{\mathcal{W}^{s}} + ||f_{P_{\geq 1}}||^{2}_{X^{s,-\q,1}} < \infty \}$$

According to the definition in (\ref{C}), the conjugate $\bar{X}^{s,b}$ is defined by:

$$\bar{X}^{s,b} =  \{f\in S'; \langle (\xi,\tau) \rangle^{s} \langle \tau+\xi^{2} \rangle^{b} \hat{f} \in L^{2}  \}$$

We can define all the other elements the same way as above by simply placing a bar on each space and operator, while replacing everywhere $|\tau-\xi^{2}|$ with $|\tau+\xi^{2}|$ and $P$ with $\bar{P}=\{(\xi,\tau): \tau+\xi^{2}=0 \}$.  

We record the following important facts:

$$ f \in X^{s,b} \Longleftrightarrow \bar{f} \in \bar{X}^{s,b} \ \ \ \mbox{and} \ \ \ (X^{0,\q})^{*}=\bar{X}^{0,-\q}$$

We owe an important remark about the use of duality. For a bilinear estimate we use the notation:

$$ \mathcal{X} \cdot \mathcal{Y} \rightarrow \mathcal{Z}$$

\noindent
which means that we seek for an estimate $||B(u,v)||_{\mathcal{Z}} \leq C ||u||_{\mathcal{X}} \cdot ||v||_{\mathcal{Y}}$. Here the constant $C$ may depend on some variables, like the frequency where the functions are localized. We use duality to claim:

$$ \mathcal{X} \cdot \mathcal{Y} \rightarrow \mathcal{Z} \Leq \mathcal{X} \cdot \mathcal{Z}^{*} \rightarrow \mathcal{Y}^{*}$$

\noindent
in the sense that the corresponding estimates are equivalent with the same constants. If $Y$ and $Z$ are not reflexive, then this is not entirely correct. Nevertheless, we apply this equivalence for functions localized in frequency, hence very smooth and belonging to the spaces we need. In this case we are allowed to use the equivalence.

For technical details we need to introduce some new localization operators. For each $i \in \N$ we define a refined lattice $\Xi^{i}= 2^{-i} \Z^{n}= \{ 2^{-i} \xi: \xi \in \Z^{n} \}$. For each $\xi \in \Xi^{i}$ we build  the corresponding $\phi^{i}_{\xi}$ to be a smooth approximation of the characteristic function of the cube centered at $\xi$ and with sizes $2^{-i-1}$. We also assume that the system $(\phi_{\xi}^{i})_{\xi \in \Xi^{i}}$ forms a partition of unity in $\R^{n}$.

For each $l \in \Z$ we can easily construct a function $\chi_{[l-\q,l+\q]}$ to be a smooth approximation of the characteristic function of the interval $[l-\q,l+\q]$ and such that the system $(\chi_{[l-\q,l+\q]})_{l \in \Z}$ form a partition of unity in $\R$. For any $\xi \in \Xi^{i}$ with $|\xi| \leq 2^{i+1}$ we consider those $l \in \Z$ with the property $|(\xi,l)| \approx 2^{i}$ and define the operators: 

$$\hat{f}_{\xi,l}= \phi_{[l-\q,l+\q]}(\tau) \phi^{i}_{\xi}(\xi) \hat{f}(\xi,\tau) $$  

  The support of $\hat{f}_{\xi, l}$ is approximately a tube centered at $(\xi,l)$ and of size $2^{-i} \times ... \times 2^{-i} \times 1$, the last one being in the $\tau$ direction. Since the distance of these tubes to $P$ will play an important role, sometimes it would be convenient if we were able to work with $(f_{\xi, \xi^{2}+l})_{\xi \in \Xi^{i}, l \in \Z}$ instead. The only problem is that it is not guaranteed that $\xi^{2} \in \Z$ for all $\xi \in \Xi^{i}$. Of course we could change the way we cut in the $\tau$ direction, but this would complicate notations even more. We choose instead to ignore that $\xi^{2}$ may not be integer, and go on and use $g_{\xi, \xi^{2}+l}$. It will be obvious from the argument that this does not affect in any way the rigorousness of the proof. The last notation we introduce is $f_{\xi, \xi^{2} \pm l}=f_{\xi, \xi^{2} + l}+f_{\xi, \xi^{2} - l}$.

For $d \leq 2^{2i-4}$ we obtain a new decomposition of $g_{i,d}$:

\beq \label{aux2}
g_{i, d} = \sum_{k =2^{-1}d}^{2 d} \sum_{\xi \in \Xi^{i}} g_{\xi, \xi^{2} \pm k}
\eeq

Notice that the $\xi$'s $\in \Xi^{i}$ involved in the above summation have $|\xi| \approx 2^{i}$.

For the part of $\hat{g}$ supported away from $P$ we come with a different decomposition: 

\beq \label{aux1}
g_{i, \geq 2^{2i-3}} = \sum_{\xi \in \Xi^{i}} \sum_{l \in I_{\xi}} g_{\xi,l}
\eeq

\noindent
where $I_{\xi}=\{l \in \Z: 2^{2i-3} \leq |l-\xi^{2}| \leq 2^{2i+2} \}$. The $\xi$'s $\in \Xi^{i}$ involved in the above summation have $|\xi| \leq 2^{i+1}$.

\section{Linear estimates} \label{les}

The linear estimates (in the case $n=2$) had been derived in \cite{be}. We briefly go over the main steps in the argument there. Nevertheless, for a full argument one should read the corresponding section in \cite{be} and notice the argument there can be easily generalized for all dimensions.

We first prove few properties which give a better insight on the structure of our spaces:

\begin{itemize}

\item $X^{s,\q,1} \subset C_{t}H^{s}_{x} \cap Z^{s} $ 

\item $ X^{s,\q,1}, Y^{s} \ \mbox{and} \ Z^s$ are stable under multiplication with $\chi_{[0,1]}$

\item $\mathcal{Y}^{s} \subset X^{s,-\q,\infty}$ and $X^{s,\q,1} \subset Y^{s}$

\item for any $d \in I_{i}$ we have $f_{i} \in Y^{s} \Ra f_{i, \leq d} \in Y^{s}$.  

\end{itemize}

Then one proves the claims for the homogeneous equation, i.e. when $f = 0$. It is well known in the literature that if $g \in H^s$ then $u$ (the solution of (\ref{I}) with $f=0$) belongs to $X^{s,b}$ for any $b$, hence it belongs to $X^{s,\q,1}$. We decompose the initial data:

$$g=\sum_{\xi \in \Xi, m \in \Z^n} g_{\xi,m}$$

\noindent
where $g_{\xi,m}$ is highly localized in a neighborhood of size $1$ around $(m,\xi)$ in phase-space. We write the solution $u$ as:

$$u=\sum_{\xi \in \Xi, m \in \Z^n} u_{\xi,m}$$

\noindent
where each $u_{\xi,m}$ solves the homogeneous equation with initial data $g_{\xi,m}$. We have the standard $L^{\infty}_{t}L^{2}_{x}$ estimates for each $u_{\xi,m}$. Then we introduce the operators $P_{j,m}(x,t,D)$, $j=1,..,n$ and $m \in \Z^{n}$:

$$P_{j,m} = x_{j}-m_{j} - i 2t D_{x_{j}} \ \ \ \mbox{with symbols} \ \ \ p_{j,m} = x_{j}-m_{j} + 2t\xi_{j}$$ 

 $i\partial_{t} - \D$ commutes with both $P_{j,m}(x,t,D)$, $j=1,..,n$ and on behalf of this we can show that $\chi_{[0,1]} u_{\xi,m}$ is highly localized in $T^{m,0}_{\xi}$. Putting everything together in order to conclude that $u \in Y^s$ is a relatively easy thing. 
 
A more delicate argument is needed for the inhomogeneous problem. Away from the paraboloid ($|\tau-\xi^2| \geq 1$) one has:

$$\hat{u}=\frac{\hat{f}}{\tau-\xi^2}$$ 

We can assume $g=0$ and decompose $f=\sum_{\xi \in \Xi, m \in \Z^n} f_{\xi,m}$ where $f_{\xi,m}$ is highly localized in $T^{m,0}_{\xi}$ (in a simplified model one can think that it is enough to take into account only the tubes with $l=0$). Along these tubes one can solve the Cauchy problem with inhomogeneity $f_{\xi,m}$ and zero initial data  to obtain the solution $u_{\xi,m}$ in $L^{\infty}_{t}L^{2}_{x}$ estimates. Using the operators $P_{j,m}(x,t,D)$ we can show that $\chi_{[0,1]} u_{\xi,m}$ is highly localized in $T^{m,0}_{\xi}$. 

With the help of the properties of our spaces, we put all the estimates together to obtain the claim in (\ref{le}).

\section{Bilinear estimates in $X^{s,\q,1}$} \label{bes}

The objective of this section is to obtain the bilinear estimates for $B(u,v)$ and $B(u,\bar{v})$ in $X^{s,\q,1}$, where $B(u,v)=\nabla{u} \cdot \nabla{v}$ . We introduce the additional bilinear form $\tilde{B}(u,v)=u \cdot v$.

The main results we claim are listed in the following theorem.

\begin{t1} \label{bil}

a) If $i \leq j$, we have the following estimates:

\beq \label{be11}
||B(u,v)||_{X^{s,-\q,1}_{k}} \leq j^2 2^{(\frac{n}{2}-s)i} 2^{(\frac{n}{2}-1+s)(k-j)}||u||_{X^{s,\q,1}_{i}} ||v||_{X^{s,\q,1}_{j}}
\eeq

b) If $2(n+2)i \leq j$, we have the following estimates:

\beq \label{be22}
||B(u,v)||_{X^{s,-\q,1}_{k, \geq 2^{k-i}}} \leq i^2 2^{(\frac{n}{2}-s)i} ||u||_{X^{s,\q}_{i}} ||v||_{X^{s,\q}_{j, \geq 2^{j-i}}}
\eeq

\end{t1}

From this result we learn that $X^{s,\q,1}$ are good spaces to estimate high-high interactions and "medium" low - high interactions. $X^{s,\q,1}$ are not good enough for the case of very low - high interactions. More rigorously, if $i \leq j$ and we can control $j$ by a small power of $2^{i}$, then the estimate (\ref{be11}) is good for our purpose. Otherwise the factor $j$ in (\ref{be11}) gets out of control. For technical reasons we decide that $2(n+2)i \leq j$ is our breakpoint threshold, i.e. we consider that if $2(n+2)i \geq j \geq i$ we keep the bilinear estimates in $X^{s,\q,1}$ and if $2(n+2)i \leq j$ then we should come up with  a different approach.

\vspace{.2in}

\subsection{Basic Estimates}
\noindent

\vspace{.1in}

We start with a simple result analyzing how two parabolas interact under convolution. We need few technical definitions. 

Throughout this section functions are defined on Fourier space (they should be thought as Fourier transforms). This is why we use the standard coordinates $(\xi,\tau)$. We say that $F$ is localized at $2^i$ iff the support of $F$ is included in the annulus $2^{i-1} \leq \langle \xi \rangle \leq 2^{i+1}$. 

For each $c \in \R$ denote by $P_{c}= \{(\xi,\tau): \tau-\xi^{2}=c\}$ and  by $\bar{P}_{c}=\{(\xi,\tau): \tau+\xi^{2}=c \}$. For simplicity $P=P_{0}$ and $\bar{P}=\bar{P}_{0}$. 

Denote by $\delta_{P_{c}}=\delta_{\tau-\xi^{2}=c}$ the standard surface measure associated to the parabola $P_{c}$. With respect to this measure, the restriction of $f$ to $P_{c}$ has norm:

$$||f||_{L^{2}(P_{c})}=\left( \int f^{2}(\xi,\xi^{2}+c) \sqrt{1+4|\xi|^{2}} d\xi  \right)^{\q}$$

Our main result in this subsections is

\begin{p1} \label{pp1}

 Let $f \in L^{2}(P^1)$ and  $g \in L^{2}(P^2)$, where $P^{1} \in \{P_{c_{1}},\bar{P}_{c_{1}}\}$ and $P^{2} \in \{P_{c_{2}},\bar{P}_{c_{2}}\}$, such that $f$ is localized at $2^{i}$ and $g$ at $2^{j}$. 
 
 a) We have:

\beq \label{ge}
||f \delta_{P^{1}} * g \delta_{P^{2}} ||_{L^{2}} \les 2^{\min{(i,j)} \frac{n}{2}} ||f||_{L^{2}(P^{1})}||g||_{L^{2}(P^{2})}
\eeq

b) If $|i-j| \leq 2$ and $k \leq \max{(i,j)}+2$ then we have:

\beq \label{ge2}
||f \delta_{P^{1}} * g \delta_{P^{2}} ||_{L^{2}(\langle \xi \rangle \approx 2^k)} \les 2^{\frac{n-1}{2}k} 2^{\frac{i}{2}} ||f||_{L^{2}(P^{1})}||g||_{L^{2}(P^{2})}
\eeq

c) If we assume $i \leq j$ and $|c_{1}|,|c_{2}| \leq 2^{i+j-5}$, then we have:

\beq \label{ge1}
||f \delta_{P^{1}} * g \delta_{P^{2}} ||_{L^{2}(|(\xi,\tau)| \approx 2^{j},|\tau-\xi^{2}| \leq d)} \les 2^{\frac{ni}{2}} (2^{-2i} d)^{\q} ||f||_{L^{2}(P^{1})}||g||_{L^{2}(P^{2})}
\eeq

\end{p1}

If we put \eqref{ge} and \eqref{ge1} together, we obtain

\begin{r2} If $i \leq j$ then

\beq \label{ge3}
||f \delta_{P^{1}} * g \delta_{P^{2}} ||_{L^{2}(|(\xi,\tau)| \approx 2^{j},|\tau-\xi^{2}| \leq d)} 
\eeq

\[
\les 2^{\frac{ni}{2}} \min{(1,(2^{-2i} d)^{\q})} ||f||_{L^{2}(P^{1})}||g||_{L^{2}(P^{2})}
\]

\end{r2}

\begin{proof}[Proof of Proposition \ref{pp1}] We make few useful simplifications. It is enough to provide the argument in the particular case when $c_1=c_2=0$ and we will explain, when needed, why this is the case. Then we can drop the $\tau$ argument for $f$ and $g$ and just think of them as functions of $\xi$; in other words rather than carrying $f(\xi,\xi^2)$ we can just use $f(\xi)$ and rewrite:

$$||f||_{L^{2}(P)}=\left( \int f^{2}(\xi) \sqrt{1+4|\xi|^{2}} d\xi  \right)^{\q}$$

a) We assume $i \leq j$ and consider two cases.

\begin{bfseries} Case 1: {\mathversion{bold}  $i \leq j-\log_{2}{\sqrt{n}}-4$} \end{bfseries} 

We perform first the argument for {\mathversion{bold} $f \delta_{P} * g \delta_{P}$}.

The strategy here is to prove an estimate of type:

\begin{equation} \label{eq1}
|(f \delta_{\tau=\xi^{2}} * g \delta_{\tau=\xi^{2}})h| \les 2^{\frac{n i}{2}} ||f||_{L^{2}(P)}||g||_{L^{2}(P)} ||h||_{L^{2}} 
\end{equation}

\noindent
for any $h \in L^{2}(\R^{3})$. We have:

$$(f \delta_{\tau=\xi^{2}} * g \delta_{\tau=\xi^{2}})h= \int h(\xi+\eta, \xi^2+\eta^2) f(\xi) g(\eta) \sqrt{1+4\xi^{2}} \sqrt{1+4\eta^{2}} d\xi d\eta $$

\noindent
where $\xi=(\xi_{1}, ..., \xi_n), \eta=(\eta_{1}, ..., \eta_{n})$. A direct use of the Cauchy-Schwarz inequality gives us:

$$|(f \delta_{\tau=\xi^{2}} * g \delta_{\tau=\xi^{2}})h| \les ||f||_{L^{2}(P)}||g||_{L^{2}(P)} \cdot $$

$$\left( \int |h|^{2}(\xi+\eta, \xi^2+\eta^2) \sqrt{1+4\xi^{2}} \sqrt{1+4\eta^{2}} d\xi d\eta\right)^{\q} \les $$

$$ 2^{\frac{i+j}{2}} ||f||_{L^{2}(P)}||g||_{L^{2}(P)} \left( \int |h|^{2}(\xi+\eta, \xi^2+\eta^2) d\xi d\eta\right)^{\q} $$

For the integral with respect to $h$ we perform the change of variables $(\xi_{1},\eta_{1},\eta_{2}, ..., \eta_n) \rightarrow (\z_{1},\z_{2}, ..., \z_{n+1})$ given by the system:

\begin{eqnarray} \label{sys3}
  \left\{
        \begin{array}{lr}
		 \xi_{1}+\eta_{1}=\z_{1}  \\
		        .....             \\
		 \xi_{n}+\eta_{n}=\z_{n}  \\
		 \xi^{2}+\eta^{2}=\z_{n+1}
	\end{array}\right.
\end{eqnarray}

This Jacobian of this transformation is $\q(\eta_{1}-\xi_{1})^{-1}$. If we were to integrate over a region where $|\eta_{1}| \geq \frac{|\eta|}{\sqrt{n}}$, then we would get $|\eta_{1} -\xi_{1}| \geq \frac{2^{j-2}}{\sqrt{n}}$ (here it is important that $i \leq j-\log_{2}{\sqrt{n}}-4$) which leads us to:

$$|(f \delta_{\tau=\xi^{2}} * g \delta_{\tau=\xi^{2}}) h| \les 2^{\frac{i}{2}} ||f||_{L^{2}(P)}||g||_{L^{2}(P)}  \left( \int ||h||_{L^{2}}^{2} d\xi_{2} ... d\xi_{n} \right)^{\q} \les$$

$$2^{\frac{ni}{2}} ||f||_{L^{2}(P)}||g||_{L^{2}(P)} ||h||_{L^{2}}$$

\noindent
the last inequality being justified by the fact that we integrate $||h||_{L^2}$ over a domain where $|\xi| \approx 2^{i}$. 

The way to fix the proof is to split $g = g_{1}+...+g_{n}$ where $g_{i}$ is localized in a region where $|\eta_{i}| \geq \frac{|\eta|}{\sqrt{n}}$. For $g_{1}$ we apply the above argument, while for $g_{i}$ we use the change of variables $(\xi_{i},\eta_{1},...,\eta_{n}) \rightarrow (\z_{1},\z_{2},...,\z_{n+1})$ given by the same system (\ref{sys3}). By adding the results we obtain for each $g_{i}$, we get (\ref{eq1}). 

In a similar way we obtain the estimates for $f \delta_{\bar{P}} * g \delta_{P}$ and $f \delta_{\bar{P}} * g \delta_{\bar{P}}$. 
 
\begin{bfseries} Case 2: {\mathversion{bold}  $ j- \log_{2}{\sqrt{n}} -3 \leq i \leq j$} \end{bfseries}

In this case we make use of the the Strichartz estimate:

$$|| \int a(\xi)  e^{i(x \cdot \xi + t \cdot \xi^{2})} d\xi  ||_{L^{4}_t L^{r}_x} \les  ||a||_{L^{2}_{\xi}}$$

\noindent
where $\frac{1}{4}+\frac{n}{2r}=\frac{n}{4}$. We would like to have a global $L^4$ estimate (i.e. $L^4_tL^4_x$) with the price of paying derivatives (powers of frequency). If $a$ is localized at $2^j$, then $b(x,t)=\int a(\xi)  e^{i(x \cdot \xi + t \cdot \xi^{2})} d\xi$ is localized in (space) frequency at $2^{j}$. Therefore we obtain:

$$||b||_{L^4_x} \les 2^{\beta j}||b||_{L^r_x}$$

\noindent
where $\beta=\frac{1}{r}-\frac{1}{4}=\frac{n}{4}-\q$. Therefore:

$$|| \int a(\xi)  e^{i(x \cdot \xi + t \cdot \xi^{2})} d\xi  ||_{L^{4}_{x,t}} \les  2^{(\frac{n}{4}-\q)j} ||a||_{L^{2}_{\xi}}$$

In our case, correcting the estimate with the measure we use on the parabola, we obtain:

$$||\mathcal{F}^{-1}(f \delta_{P})||_{L^{4}} \leq 2^{\frac{n i}{4}} ||f||_{L^{2}(P)} \ \ \mbox{and} \ \ ||\mathcal{F}^{-1}(g \delta_{P})||_{L^{4}} \leq 2^{\frac{n j}{4}} ||g||_{L^{2}(P)}$$

We can conclude with: 

$$||f \delta_{P}* g \delta_{P}||_{L^{2}} = ||\mathcal{F}^{-1}(f \delta_{P}) \cdot \mathcal{F}^{-1}(g \delta_{P})||_{L^{2}} \les$$

$$||\mathcal{F}^{-1}(f \delta_{P})||_{L^{4}}||\mathcal{F}^{-1}(g \delta_{P})||_{L^{4}} \les 2^{\frac{nj}{2}} ||f||_{L^{2}(P)} ||g||_{L^{2}(P)}$$

Since the Strichartz estimate is valid for $\bar{P}$ too, we can extend the argument to the other combinations, i.e.  $f \delta_{\bar{P}} * g \delta_{P}$ and $f \delta_{\bar{P}} * g \delta_{\bar{P}}$. 

\vspace{.1in}

b) This result is a refinement of the one in a) for the case when the interaction needs to be measured on smaller sets. We say smaller, since there is nothing to prove if, let's say, $|k-j| \leq 10$. Therefore we assume that $k \leq j-10$ for the rest of this part of the argument. We start the same way as in part a), i.e., we want to estimate:

$$(f \delta_{\tau=\xi^{2}} * g \delta_{\tau=\xi^{2}})h= \int h(\xi+\eta, \xi^2+\eta^2) f(\xi) g(\eta) \sqrt{1+4\xi^{2}} \sqrt{1+4\eta^{2}} d\xi d\eta $$

\noindent
for $h \in L^{2}$ supported at $2^k$. We decompose $\R^n$ in disjoint cubes of size $2^{k+1}$:

$$\R^{n}=\cup_{l \in \Z^n} C_{l}$$

\noindent
where $C_l$ is a cube of size $2^{k+1}$, centered at $2^{k+1} l$ and with the sides parallel to the standard coordinate axes. We denote by $f_l$ the part of $f$ supported in $C_l$ and similarly for $g_l$. We observe that for each $l$ the set $B_{l}=\{l' \in \Z^n: (C_l + C_{l'}) \cap C_0 \ne \emptyset\}$ has cardinality at most $\approx 2^n+1$. Also for each $l'$ there are at most $\approx 2^n+1$ values of $l$ such that $l' \in B_{l}$. 

We have:

$$(f \delta_{\tau=\xi^{2}} * g \delta_{\tau=\xi^{2}})h= \sum_{l \in \Z^n} \sum_{l' \in B_l} \int h(\xi+\eta, \xi^2+\eta^2) f_l(\xi)  g_{l'}(\eta) \sqrt{1+4\xi^{2}} \sqrt{1+4\eta^{2}} d\xi d\eta $$

For fixed $l$ and $l'$ we run the same argument as in part a), Case 1,  up to the point:

$$|(f_l \delta_{\tau=\xi^{2}} * g_{l'} \delta_{\tau=\xi^{2}}) h| \les 2^{\frac{i}{2}} ||f_l||_{L^{2}(P)}||g_{l'}||_{L^{2}(P)}  \left( \int ||h||_{L^{2}}^{2} d\xi_{2} ... d\xi_{n} \right)^{\q} $$

Since we deal with the case when $k \leq j-10$, it follows that if $\xi \in C_l$, $\eta \in C_{l'}$ and $\xi +\eta \in C_{0}$, then $|\eta - \xi| \geq 2^{j-5}$, hence we need to run only the argument in Case 1 from part a). In the domain of integration above we have $\Delta \xi_{p} \approx 2^{k}$ for all $2 \leq p \leq n$; therefore we obtain:

$$|(f_l \delta_{\tau=\xi^{2}} * g_{l'} \delta_{\tau=\xi^{2}}) h| \les 2^{\frac{i}{2}} 2^{\frac{n-1}{2}k} ||f_l||_{L^{2}(P)}||g_{l'}||_{L^{2}(P)}  ||h||_{L^{2}} $$

Summing up with respect to $l$ and $l' \in B_l$  and applying Cauchy-Schwarz gives us:

$$|(f \delta_{\tau=\xi^{2}} * g \delta_{\tau=\xi^{2}}) h| \les 2^{\frac{i}{2}} 2^{\frac{n-1}{2}k} ||f||_{L^{2}(P)}||g||_{L^{2}(P)}  ||h||_{L^{2}} $$

In obtaining this estimate we used the fact the observations made at the beginning of the argument about the sets $B_{l}$.

\vspace{.1in}

c) It is enough to prove this result under the hypothesis that $d \leq 2^{i+j-4}$ since otherwise, the result in (\ref{ge}) is stronger. Without losing generality we can assume $c_{1}=c_{2}=0$. One could easily adapt the argument bellow to the general case when $|c_1|+|c_2| \leq 2^{i+j-4}$. 

{\mathversion{bold} $f \delta_{P} * g \delta_{P}$}

We test $||f \delta_{\tau=\xi^{2}} * g \delta_{\tau=\xi^{2}}||_{L^2(|\tau-\xi^2| \leq d)}$  by estimating $|(f \delta_{\tau=\xi^{2}} * g \delta_{\tau=\xi^{2}})h|$ for any $h \in L^{2}$ supported in the region $|\tau-\xi^{2}| \leq d$. For any such $h$ we have:

$$(f \delta_{\tau=\xi^{2}} * g \delta_{\tau=\xi^{2}})h=$$

$$\int f(\xi) g(\eta) h(\xi+\eta, \xi^{2}+\eta^{2}) \sqrt{1+4\xi^{2}}  \sqrt{1+4\eta^{2}} d\xi d\eta $$

Since $h$ is supported in a region $|\tau-\xi^{2}| \leq d$ we need the following condition on the variables inside the integral: $|(\xi+\eta)^{2}-(\xi^{2}+\eta^{2})| \leq d$ or $2|\xi||\eta|\cos{\t} \leq d$ where $\t$ is the angle between $\xi$ and $\eta$. 
Hence $\cos{\t} \leq 2^{-i-j}d$ which implies $|\t-\frac{\pi}{2}| \leq 2^{-i-j}d$. This suggests decomposing $\R^n$ in conical sets which have angular dimension $\a^{n}$ where $\a=2^{-i-j}d$. To formalize this we write:

$$\R^n=\cup_{l \in S_{\a}} A_{l}$$

\noindent
where $S_{\a}$ contains $\approx \a^{1-n}$ indexes and, for any $l \in S_{\a}$, any two vectors in $A_{l}$ make an angle of at most $\a$. Correspondingly we split:

$$f=\sum_{l \in S_{\a}} f_{l}  \ \ \ \mbox{and} \ \ \ g=\sum_{l \in S_{\a}} g_{l}$$

\noindent
such that $f_{l}$ is the part of $f$ localized in $A_{l}$ and similarly for $g$.

For indexes in $S_{\a}$ we define the following relation:

$$l \perp l' \Leq \exists \xi \in A_{l}, \eta \in A_{l'} \ \mbox{such that} \ |\angle (\xi,\eta) - \frac{\pi}{2}| \leq \a$$

\noindent
where $\angle (\xi,\eta)$ denotes the angle between $\xi$ and $\eta$. 

We notice the following property: for each $l$ there are about $\a^{2-n}$ $l'$'s such that $l \perp l'$. To understand better this observation, one could think of $\cup_{l' \perp l} A_{l'}$ as a hyperplane of codimension $1$ since it contains "essentially" the vectors which are orthogonal on the ones in $A_l$; this way it loses one degree of freedom in the set $S$ and then it contains $\approx \a^{2-n}$ elements.

If $\xi \in A_{l}$ and $ \eta \in A_{l'}$ and we want them to belong to the domain of integration above we need to impose $l \perp l'$. Then we have:

\beq \label{ba1}
(f \delta_{\tau=\xi^{2}} * g \delta_{\tau=\xi^{2}})h=
\eeq

$$\sum_{l \in A_{\a}} \int f_{l}(\xi) \sum_{l' \perp l} g_{l'}(\eta) h_{l,l'}(\xi+\eta, \xi^{2}+\eta^{2}) \sqrt{1+4\xi^{2}}  \sqrt{1+4\eta^{2}} d\xi d\eta $$

For fixed $l$ we notice that the supports of $\{\xi+\eta: \xi \in A_l, \eta \in A_{l'} \}$ are "essentially" disjoint with respect to $l' \perp l$. This is exactly the projection of the support of $h_{l,l'}$ onto the $\xi$ plane.

Now that we have a sharp angular localization, we complete it with a norm localization which should be consistent with the angular one:

$$f_{l}=\sum_{m} f_{l,m} \ \ \ \mbox{and} \ \ \ g_{l'}=\sum_{n} g_{l',n}$$

We describe bellow the construction of this decompositions.

We pick two vectors $v_{l} \in A_{l}$ and $v_{l'} \in A_{l'}$ such that $v_{l} \perp v_{l'}$. In what follows we will indicate the sizes of a  parallelepiped with respect to $v_{l'} \times v_{l} \times ....$, to mean that the first one is the size in the direction of $v_{l'}$, the second one is the size with respect to $v_l$ and the others are with respect to the orthogonal directions to $v_{l'}$ and $v_{l}$.

$f_{l,m}$ is the part of $f_{l}$ localized in the set  $A^{i}_{l,m}=\{\xi \in A_{l}:|\xi| \in [2^{-i}d(m-\q),2^{-i}d(m+\q)] \}$ and notice that this is consistent with the arc length size localization of $g_{l'}$ (which is $2^{-i}d$). The support of $A^{i}_{l,m}$ is 
approximately a parallelepiped of sizes $2^{-j}d \times 2^{-i}d \times (2^{-j}d)^{n-2}$ with respect to $v_{l'} \times v_{l} \times ...$; last dimension should be understood as $2^{-j}d \times ... \times 2^{-j}d$ ($n-2$ times).

For $g_{l'}$ we want to do something similar: we would like to localize $|\eta|$ in intervals of size $2^{-j}d$. The only problem we encounter is that if $i << j$ and $d$ small we may see the curvature of the sphere and then the support of $g_{l',n}$ cannot be approximated by a rectangle.   

In order to fix this we chose $g_{l',n}$ to be the part of $g_{l}$ localized in  $A^{j}_{l',n}=\{\eta \in A_{l^{\perp}}: \eta \cdot v_{l'} \in [2^{-j}d(n-\q),2^{-j}d(n+\q)] \}$. The support of $A^{j}_{l',n}$ is approximately a parallelepiped of sizes $2^{-j}d \times 2^{-i}d \times (2^{-i}d)^{n-2}$.

The crucial property is that the sum sets of the supports, namely $A^{i}_{l,m}+A^{j}_{n,l^{\perp}}=\{\xi+\eta:\xi \in A^{i}_{l,m} \ \mbox{and} \ \eta \in  A^{j}_{l',n}\}$ are disjoint with respect to the pair $(m,n)$. This is because the sum set $A^{i}_{l,m}+A^{j}_{n,l'}$ is approximately a parallelepiped of sizes $2^{-j}d \times 2^{-i}d \times (2^{-i}d)^{n-2}$ and whose center has coordinates $(2^{-j}dn,2^{-i}dm, ....)$ with respect to the base $v_{l'} \times v_{l} \times ...$. We denote by $h_{l,l',m,n}$ the part of $h$ which is supported in this set (more precisely the projection of the support on the $\xi$ space should be supported there). Hence we can write:

$$(f \delta_{\tau=\xi^{2}} * g \delta_{\tau=\xi^{2}})h=$$

$$\sum_{l \in S_{\a}} \sum_{l' \perp l} \sum_{m} \sum_{n} \int f_{l,m}(\xi) g_{l',n}(\eta) h_{l,l',m,n}(\xi+\eta, \xi^{2}+\eta^{2}) \sqrt{1+4\xi^{2}}  \sqrt{1+4\eta^{2}} d\xi d\eta $$

 For fixed $l$ and $l'$, we can estimate, via Cauchy-Schwarz:

$$|(f_{l} \delta_{\tau=\xi^{2}} * g_{l'} \delta_{\tau=\xi^{2}})h| \les  \left( \sum_{m} ||f_{l,m}||^{2}_{L^{2}(P)} \right)^{\q} \left( \sum_{n} ||g_{l',n}||^{2}_{L^{2}(P)} \right)^{\q} \cdot$$

$$\left( \sum_{m,n} \int h^{2}_{l,l',m,n}(\xi+\eta,\xi^{2}+\eta^{2}) \sqrt{1+4\xi^{2}}  \sqrt{1+4\eta^{2}} d\xi d\eta \right)^{\q}$$

We can perform a change of coordinates whose Jacobian is $1$ such that $\xi_1$ and $\eta_1$ become the coordinates in the direction of $v_{l'}$ and $\xi_2$ and $\eta_2$ become the coordinates in the direction of $v_{l}$. Once we achieve that, we notice that $|\eta_1-\xi_1| \approx 2^{j}$ and then we can perform the change of coordinates we introduced in (\ref{sys3}) and estimate:

$$\int h^{2}_{l,l',m,n}(\xi+\eta,\xi^{2}+\eta^{2}) \sqrt{1+4\xi^{2}}  \sqrt{1+4\eta^{2}} d\xi d\eta \les 2^i \int ||h_{l,l',m,n}||_{L^{2}}^{2} d\xi_{2} ... d\xi_{n} \les $$ 

$$2^i 2^{j} \a (2^i \a)^{n-2}  ||h_{l,l',m,n}||_{L^{2}}^{2} $$

In the last estimate we used that $\Delta \xi_2 \approx 2^j \a$ and $\Delta \xi_k \approx 2^i \a$ for each $3 \leq k \leq n$ (recall the dimensions of the $A^{i}_{l,m}$). Since the supports of $h_{l,l',m,n}$ are disjoint with respect to the pair $(m,n)$, we can conclude the above computation with:

$$|(f_{l} \delta_{\tau=\xi^{2}} * g_{l'} \delta_{\tau=\xi^{2}})h_{l,l'}| \les  2^{\frac{i}{2}} (2^{j}\a)^{\q}(2^i \a)^{\frac{n-2}{2}}||f_{l}||_{L^{2}(P)}  ||g_{l'}||_{L^{2}(P)} ||h_{l,l'}||_{L^{2}}$$

Recalling the fact that the supports of $h_{l,l'}$ are disjoint with respect to $l' \perp l$, we can perform the summation with respect to $l' \perp l$ ($l$ is fixed):

$$|\sum_{l' \perp l} (f_{l} \delta_{\tau=\xi^{2}} * g_{l'} \delta_{\tau=\xi^{2}})h_{l,l'}| \les 2^{\frac{i}{2}} (2^{j}\a)^{\q}(2^i \a)^{\frac{n-2}{2}}||f_{l}||_{L^{2}(P)}  ||\sum_{l' \perp l} g_{l'}||_{L^{2}(P)} ||\sum_{l' \perp l} h_{l,l'}||_{L^{2}} $$

$$ \les 2^{\frac{i}{2}}(2^{j}\a)^{\q}(2^i \a)^{\frac{n-2}{2}}||f_{l}||_{L^{2}(P)}  ||\sum_{l' \perp l} g_{l'}||_{L^{2}(P)} ||h||_{L^{2}}$$

In the end we perform the summation with respect to $l$:

$$|(f \delta_{\tau=\xi^{2}} * g \delta_{\tau=\xi^{2}})h| \les \sum_{l} 2^{\frac{i}{2}}(2^{j}\a)^{\q}(2^i \a)^{\frac{n-2}{2}} ||f_{l}||_{L^{2}(P)}  ||\sum_{l' \perp l} g_{l'}||_{L^{2}(P)} ||h||_{L^{2}} \les$$

$$ 2^{\frac{i}{2}}(2^{j}\a)^{\q}(2^i \a)^{\frac{n-2}{2}}\left( \sum_{l} ||f_{l}||^{2}_{L^{2}(P)} \right)^{\q} \left( \sum_{l} ||\sum_{l' \perp l}g_{l'}||^{2}_{L^{2}(P)} \right)^{\q} ||h||_{L^{2}} \les $$

$$ 2^{\frac{n-1}{2}i} (2^{j}\a)^{\q} ||f||_{L^{2}(P)} ||g||_{L^{2}(P)} ||h||_{L^{2}}$$

In the last estimate we have used the fact (mentioned at the beginning) that for every $l'$ there are about $\a^{2-n}$ values of $l$ such that $l \perp l'$. 

Since this holds true for any $h \in L^{2}$ supported in $|\tau-\xi^{2}| \leq d$, we can conclude with the result of the Proposition.

{\mathversion{bold} $f \delta_{\bar{P}} * g \delta_{P}$}

We start the argument in a similar way. We test the convolution against a $h \in L^{2}$ supported in $|\tau-\xi^{2}| \leq d$:

$$(f \delta_{\tau=-\xi^{2}} * g \delta_{\tau=\xi^{2}})h=$$

$$\int f(\xi) g(\eta) h(\xi+\eta, -\xi^{2}+\eta^{2}) \sqrt{1+4\xi^{2}}  \sqrt{1+4\eta^{2}} d\xi d\eta $$

Since $h$ is supported in a region $|\tau-\xi^{2}| \leq d$ we need the following condition on the variables inside the integral: $|(\xi+\eta)^{2}-(-\xi^{2}+\eta^{2})| \leq d$ or $2|\xi||\eta+\xi|\cos{\t} \leq d$ where $\t$ is the angle between $\xi$ and $\eta+\xi$. Hence $\cos{\t} \leq 2^{-i-j}d$ which implies $|\t-\frac{\pi}{2}| \leq 2^{-i-j}d$. This suggests decomposing:

$$f=\sum_{l \in S_{\a}} f_{l}  \ \ \ \mbox{and} \ \ \ h=\sum_{l \in S_{\a}} h_{l}$$

\noindent
such that $f_{l}$ is as before and $h_{l}$ is the part of $h$ whose support, when projected on the $\xi$ plane, is included in $A_{l}$. Then the argument continues as before with $h$ taking the place of $g$ and vice-versa. 

{\mathversion{bold} $f \delta_{P} * g \delta_{\bar{P}}$}

If $i \leq j-5$ then the convolution is localized in region with $\tau \leq 0$, hence outside the region with $|\tau-\xi^{2}| \leq d$. If $j-5 \leq i \leq j$ then this is similar to the case {\mathversion{bold} $f \delta_{\bar{P}} * g \delta_{P}$}.

\end{proof}

\vspace{.2in}

\subsection{Bilinear estimates on dyadic regions}
\noindent
\vspace{.1in}

A standard way of writing down each case looks like:

\vspace{.1in}

{\mathversion{bold} $ X_{i, d_{1}}^{0,\q} \cdot X_{j, d_{2}}^{0,\q} \rightarrow  X_{j, d_{3}}^{0,-\q}$  }

\vspace{.1in}

This means that for $u \in X_{i, d_{1}}^{0,\q}$ and $v \in X_{j, d_{2}}^{0,\q}$ we estimate the part of $B(u,v)$ (or $\tilde{B}(u,v)$) whose Fourier transform is supported in $A_{j, d_{3}}$. Formally we estimate $\mathcal{F}^{-1}(\chi_{A_{j,d_{3}}}\mathcal{F}(B(u,v)))$. This is going to be the only kind of ``abuse'' in notation which we make throughout the paper, i.e. considering $||B(u,v)||_{X_{j,d_{3}}^{s,\q}}$ even if $\mathcal{F}(B(u,v))$ is not supported in $A_{j,d_{3}}$. We choose to do this so that we do not have to relocalize every time in $A_{j,d_{3}}$.

 Sometimes we prove estimates via duality or conjugation by using this simple properties: 

$$X \cdot Y \rightarrow Z \Longleftrightarrow X \cdot (Z)^{*} \rightarrow (Y)^{*} \ \ \mbox{and} \ \  X \cdot Y \rightarrow Z \Longleftrightarrow \bar{X} \cdot \bar{Y} \rightarrow \bar{Z}$$

\begin{p1} \label{u1}

Assume $0 \leq i \leq j$ and $|k-j| \leq 4$. Then we have the estimates:

\beq \label{b1}
||\tilde{B}(u,v)||_{X^{0,-\q}_{k,d_{3}}} \les 2^{\frac{n-1}{2}i} 2^{-\frac{j}{2}} \min{(2^{-i}, d_{2}^{-\q},   d_{3}^{-\q})} ||u||_{X^{0,\q}_{i,d_{1}}} ||v||_{X^{0,\q}_{j,d_{2}}}
\eeq

If $\max{(d_{2},d_{3})} \geq 2^{i+j+6}$, then we have a nontrivial estimate only if $ 4^{-1} \leq d_{2}d_{3}^{-1} \leq 4$ and the estimate is improved to:

\beq \label{b2}
||\tilde{B}(u,v)||_{X^{0,-\q}_{k,d_{3}}} \les  2^{\frac{n}{2} i} d_{3}^{-1} ||u||_{X^{0,\q}_{i,d_{1}}} ||v||_{X^{0,\q}_{j,d_{2}}}
\eeq

\noindent
without any restriction over $d_{1}$. 

Assume $|i-j| \leq 2$ and $k \leq j-1$. Then we have the estimates:

\beq \label{b10}
||\tilde{B}(u,v)||_{X^{0,-\q}_{k,d_{3}}} \les 2^{\frac{n-1}{2}k} 2^{-\frac{j}{2}} \min{(2^{-k},d_{1}^{-\q}, d_{2}^{-\q})} ||u||_{X^{0,\q}_{i,d_{1}}} ||v||_{X^{0,\q}_{j,d_{2}}}
\eeq

Except (\ref{b2}), all of the above estimates hold true if $\tilde{B}(u,v)$ is replaced by $\tilde{B}(\bar{u},v)$ or $\tilde{B}(u,\bar{v})$. (\ref{b2}) holds true, with the same restrictions, if $\tilde{B}(u,v)$ is replaced by $\tilde{B}(\bar{u},v)$.

\end{p1}

In order to prove the estimates (\ref{A}) and (\ref{M}) we need the following results which are straightforward consequences from the corresponding ones in Proposition \ref{u1}.

\begin{p1} \label{u9}

Assume $0 \leq i \leq j+2$ and $|k-j| \leq 4$. Then we have the estimates:

\beq \label{b99}
||\tilde{B}(u,v)||_{X^{0,\q}_{k,d_{3}}} \les 2^{\frac{n-1}{2}i} 2^{-\frac{j}{2}} d_{3}^{\q} ||u||_{X^{0,\q}_{i,d_{1}}} ||v||_{X^{0,\q}_{j,d_{2}}}
\eeq

If $\max{(d_{2},d_{3})} \geq 2^{i+j+6}$, then we have a nontrivial estimate in (\ref{b99}) only if $ 4^{-1} \leq d_{2}d_{3}^{-1} \leq 4$ and the estimate is improved to:

\beq \label{b96}
||\tilde{B}(u,v)||_{X^{0,\q}_{k,d_{3}}} \les  2^{\frac{n}{2}i} ||u||_{X^{0,\q}_{i,d_{1}}} ||v||_{X^{0,\q}_{j,d_{2}}}
\eeq

If $|i-j| \leq 2$ we have the estimates:

\beq \label{b95}
||\tilde{B}(u,v)||_{X^{0,\q}_{k,d_{3}}} \les 2^{(\frac{n}{2}-2)k} d_{3} ||u||_{X^{0,\q}_{i,d_{1}}} ||v||_{X^{0,\q}_{j,d_{2}}}
\eeq

\beq \label{b195}
||\tilde{B}(u,v)||_{X^{0,-\q}_{k,d_{3}}} \les 2^{(\frac{n}{2}-2)k} ||u||_{X^{0,\q}_{i,d_{1}}} ||v||_{X^{0,\q}_{j,d_{2}}}
\eeq

(\ref{b99}), (\ref{b95}) and (\ref{b195}) hold true if $B(u,v)$ is replaced by $B(\bar{u},v)$ or $B(u,\bar{v})$, (\ref{b96}) holds true if we replace $\tilde{B}(u,v)$ by $\tilde{B}(\bar{u},v)$.
\end{p1}

We introduce the following operator: $\hat{B}(u,v)=\mathcal{F}^{-1}({|\hat{u}| * |\hat{v}|})$. $\hat{B}$ is similar to $\tilde{B}$ just that it takes absolute values on the Fourier side. On the other hand, one can easily see that in the above estimates we pass immediately to absolute values on the Fourier side, hence every estimate from above holds true with $\tilde{B}$ replaced by $\hat{B}$. The reason we need $\hat{B}$ is the following  

\begin{p1} \label{u0} If $u \in L^{2}_{i,d_{1}}$ and $v \in L^{2}_{j,d_{2}}$, we have the estimate:

\beq \label{r1}
||B(u,v)||_{L^{2}_{k,d_{3}}} \les \sup_{\xi, \eta} |\xi \cdot \eta| \cdot ||\hat{B}(u,v)||_{L^{2}}
\eeq

\noindent
where the $\sup$ is taken over the pairs $(\xi,\tau_1) \in A_{i,d_{1}}$, $(\eta,\tau_2) \in A_{j,d_{2}}$ such that $(\xi+\eta,\tau_1+\tau_2) \in A_{k,d_{3}}$. 

\end{p1}

\begin{proof}

By taking the Fourier transform we have:

$$\widehat{B(u,v)}(\z,\tau_{3})=-\int \xi \cdot \eta \hat{u}(\xi,\tau_{1}) \hat{v}(\eta,\tau_{2}) d \xi d \tau_1$$

\noindent
where the domain of the integral is restricted to where $\xi+\eta=\z$ and $\tau_1 +\tau_2=\tau_3$. Then the conclusion is obvious.

\end{proof}

We could not replace in the above estimate $\hat{B}(u,v)$ with $\tilde{B}(u,v)$. What we are allowed instead is to estimate $\hat{B}(u,v)$ by the bounds obtained in Proposition \ref{u1} for $\tilde{B}(u,v)$.  

We have the identity  $\tau_1+\tau_2-(\xi +\eta)^2=\tau_1 - \xi^2 + \tau_2 -\eta^2 -2 \xi \cdot \eta$. Therefore, if $u \in L^{2}_{i,d_{1}}$, $v \in L^{2}_{j,d_{2}}$ and we want to estimate $B(u,v)$ in $L^{2}_{k,d_{3}}$ then we have:

\beq \label{res1}
\sup_{\xi, \eta} |\xi \cdot \eta| \les \max{(d_1,d_2,d_3)}
\eeq

Now we can adapt the results in Proposition \ref{u1} for the operator $B$:

\begin{p1} \label{u4}

Assume $0 \leq i \leq j$ and $|k-j| \leq 4$. Then we have the estimates:

\beq \label{b100}
||B(u,v)||_{X^{0,-\q}_{j,d_{3}}} \les 2^{\frac{n-1}{2}i} 2^{-\frac{j}{2}} \max{(2^{-i}d_1, d_{2}^{\q},   d_{3}^{\q})}  ||u||_{X^{0,\q}_{i,d_{1}}} ||v||_{X^{0,\q}_{j,d_{2}}}
\eeq

If $\max{(d_{2},d_{3})} \geq 2^{i+j+6}$, then we have a nontrivial estimate only if $ 4^{-1} \leq d_{2}d_{3}^{-1} \leq 4$ and the estimate is improved to:

\beq \label{b102}
||B(u,v)||_{X^{0,-\q}_{j,d_{3}}} \les  2^{\frac{n}{2}i} ||u||_{X^{0,\q}_{i,d_{1}}} ||v||_{X^{0,\q}_{j,d_{2}}}
\eeq

\noindent
without any restriction over $d_{1}$. 

If $|i-j| \leq 2$ and $k \leq j-1$, then we have the estimates:

\beq \label{b103}
||\tilde{B}(u,v)||_{X^{0,-\q}_{k,d_{3}}} \les 2^{\frac{n-1}{2}k} 2^{-\frac{j}{2}} \min{(2^{-k}d_3,d_{1}^{\q}, d_{2}^{\q})} ||u||_{X^{0,\q}_{i,d_{1}}} ||v||_{X^{0,\q}_{j,d_{2}}}
\eeq

\end{p1}

\begin{proof}[Proof of Proposition \ref{u1}]

We make some commentaries about the statement above. If $i \leq j-2$, then the result is localized at frequency $\approx 2^{j}$. There is something to estimate only if $k=j, j -1$. 

It is only when $i=j-1,j$ that we have parts of the result at lower frequencies and then we have to provide estimates for all $k \leq j+1$.

We deal first with the case when we measure the outcome at the high frequency and at the end we deal with the case when we have $i=j-1,j$ and we have to measure the outcome at lower frequencies. 

We need  to transform the estimates on paraboloids in estimates on dyadic pieces. If we localize in a region where $ \langle \xi \rangle \approx 2^{k}$, the paraboloids $P_{c}$ make an angle of $\approx 2^{-k}$ with the $\tau$ axis, so we have the following relation between measures:

\beq \label{meas}
d\xi d\tau \approx 2^{-k} dP_{c} dc
\eeq

If $d \leq 2^{2i-4}$ then for $(\xi,\tau) \in A_{i,d}$ we have $|\xi| \approx 2^{i}$. Therefore for $d \leq 2^{2i-4}$:

\beq \label{e100}
||u||^{2}_{X^{0,\pm \q}_{i,d}} \approx (1+d)^{\pm \q} \int_{b=\frac{d}{2}}^{2d} ||\hat{u}||^{2}_{L^{2}(P_{b})} 2^{-i} db 
\eeq

\vspace{.1 in}

\begin{bfseries} Low - High interactions with output at high frequencies  \end{bfseries} 

In order to be more suggestive about the kind of estimates we provide in this part of the proof, we choose to list them only for $k=j$. One would easily notice that we can replace the right hand side space with $X_{k, d_{3}}^{0,-\q}$ for any $|k-j| \leq 4$ in the estimates bellow. 

\vspace{.1 in}

{\mathversion{bold} $ X_{i, d_{1}}^{0,\q} \cdot  X_{j, d_{2}}^{0,\q} \rightarrow  X_{j, d_{3}}^{0,-\q} \ \ \mbox{and} \ \ \bar{X}_{i, d_{1}}^{0,\q} \cdot X_{j, d_{2}}^{0,\q} \rightarrow X_{j, d_{3}}^{0,-\q}$}

\vspace{.1 in}

\begin{bfseries} Case 1: $d_{1} \leq 2^{2i-4}$ \end{bfseries}

\begin{bfseries} Subcase 1.1: $d_{2} \leq 2^{2j-4}$ \end{bfseries}

We use \eqref{meas} and then apply (\ref{ge3}) to estimate:

$$||\hat{u} * \hat{v}||_{L^{2}_{j,d_3}} \leq \int_{I_{1}} \int_{I_{2}} ||\hat{u} \delta_{P_{b_{1}}} * \hat{v} \delta_{P_{b_{2}}} ||_{L^{2}_{j,d_3}} 2^{-i-j} db_{1} db_{2} \leq $$

$$ \int_{I_{1}} \int_{I_{2}} 2^{\frac{n-2}{2}i} 2^{-j} \min{(1,(2^{-2i}d_3)^\q)} ||\hat{u}||_{L^{2}(P_{b_{1}})} ||\hat{v}||_{L^{2}(P_{b_{2}})} db_{1} db_{2} \leq $$

$$  2^{\frac{n-1}{2}i} 2^{\frac{-j}{2}} \min{(1,(2^{-2i}d_3)^\q)}\left( \int_{I_{1}} (1+b_{1})^{-1} db_{1} \right)^{\q} ||u||_{X^{0,\q}_{i,d_{1}}} \left( \int_{I_{2}} (1+b_{2} )^{-1} db_{2} \right)^{\q} ||v||_{X^{0,\q}_{j,d_{2}}} \approx $$

$$ 2^{\frac{n-1}{2}i} 2^{\frac{-j}{2}} \min{(1,(2^{-2i}d_3)^\q)}||u||_{X^{0,\q}_{i,d_{1}}} ||v||_{X^{0,\q}_{j,d_{2}}}$$

Here we used the fact that $I_{1} \approx [\frac{d_{1}}{2},2d_{1}]$ which gives us $\int (1+b_{1})^{-1} db_{1} \approx 1 $; same thing applies  for the integral with respect to $b_{2}$. Then

$$||\tilde{B}(u,v)||_{X^{0,-\q}_{j,d_{3}}} \leq d_{3}^{-\q} ||\hat{u} * \hat{v}||_{L^{2}_{j,d_3}} \leq  2^{\frac{n-1}{2}i} 2^{\frac{-j}{2}} \min{(2^{-i},d_{3}^{-\q})} ||u||_{X^{0,\q}_{i,d_{1}}} ||v||_{X^{0,\q}_{j,d_{2}}}$$ 

 Since the results in (\ref{ge}) and (\ref{ge1}) allow us to replace $\delta_{P_c}$ with $\delta_{\bar{P}_c}$ we obtain in a similar way the estimates for $\bar{X}_{i, d_{1}}^{0,\q} \cdot X_{j, d_{2}}^{0,\q} \rightarrow X_{j, d_{3}}^{0,-\q}$.

\begin{bfseries} Subcase 1.2: $d_{3} \leq 2^{2j-4}$ \end{bfseries}

This estimate for this case can be deduced by duality from the estimate:

$$ X_{i, d_{1}}^{0,\q} \cdot  \bar{X}_{j, d_{3}}^{0,\q} \rightarrow  \bar{X}_{j, d_{2}}^{0,-\q} \Leftrightarrow  \bar{X}_{i, d_{1}}^{0,\q} \cdot  X_{j, d_{3}}^{0,\q} \rightarrow  X_{j, d_{2}}^{0,-\q} $$

\begin{bfseries} Subcase 1.3: $\max{(d_{2},d_{3})} \geq 2^{i+j+6}$  \end{bfseries}

In this case we prefer to prove directly \eqref{b2} since it gives us a stronger result then \eqref{b1}.

Let $(\xi_{1},\tau_{1}) \in A_{i,d_{1}}$ and $(\xi_{2},\tau_{2}) \in A_{j,d_{2}}$. Their interaction, under convolution is $(\xi_1+\xi_2,\tau_{1}+\tau_{2})$ and satisfies $|\tau_1+\tau_2-(\xi_{1}+\xi_2)^{2}|=|\tau_1-\xi_1^{2}+\tau_{2}-\xi_{2}^{2}-2\xi_{1}\xi_{2}| \approx |\tau_{2}-\xi_{2}^{2}|$, since $|\tau_1-\xi_{1}^2| \leq 2^{2i+2} \leq 2^{i+j+3}$ and $|2 \xi_1 \xi_2| \leq 2^{i+j+3}$. Therefore the support of $\hat{u} * \hat{v}$ is included in $\cup_{d_{3}=2^{-2}d_{2}}^{2^{2}d_{2}} A_{j,d_2}$.   

In this case we use a much simpler argument, which, for reference, we call the $L^{1} * L^{2} \rightarrow L^{2}$ argument. It goes as follows:

$$||\hat{u}||_{L^{1}} \leq 2^{\frac{n}{2}i} d_{1}^{\q} ||\hat{u}||_{L^{2}} \leq  2^{\frac{n}{2}i} ||u||_{X_{i,d_{1}}^{0,\q}} $$

Then we continue with:

$$||\tilde{B}(u,v)||_{X^{0,-\q}_{j,d_{3}}} \approx  d_{3}^{-\q} ||\hat{u} * \hat{v} ||_{L^{2}(A_{j,d_{3}})} \leq $$

$$ d_{3}^{-\q} ||\hat{u}||_{L^{1}} \cdot ||\hat{v}||_{L^{2}} \leq 2^{\frac{n}{2}i} d_{3}^{-1} ||u||_{X^{0,\q}_{j,d_{1}}} ||v||_{X^{0,\q}_{j,d_{2}}}$$

A similar approach gives us the estimate for $\bar{X}_{i, d_{1}}^{0,\q} \cdot X_{j, d_{2}}^{0,\q} \rightarrow X_{j, d_{3}}^{0,-\q}$.

\begin{bfseries} Subcase 1.4: $\min{(d_{2},d_{3})} \geq 2^{2j-3}$  \end{bfseries}

If $i$ is close to $j$ then there is a gap between the ranges in the subcases 1.1 and 1.2, on one hand, and the subcase 1.3, on the other hand. If $d_3 \geq d_2$, we apply the $L^{1} * L^{2} \rightarrow L^{2}$:

$$||\tilde{B}(u,v)||_{X^{0,-\q}_{j,d_{3}}} \approx  d_{3}^{-\q} ||\hat{u} * \hat{v} ||_{L^{2}(A_{j,d_{3}})} \leq $$

$$ d_{3}^{-\q} ||\hat{u}||_{L^{1}} \cdot ||\hat{v}||_{L^{2}} \leq 2^{\frac{n}{2}i} d_{3}^{-\q} 2^{-j} ||u||_{X^{0,\q}_{j,d_{1}}} ||v||_{X^{0,\q}_{j,d_{2}}}$$

If $d_3 \leq d_2$ we derive the estimate via duality from:

$$ X_{i, d_{1}}^{0,\q} \cdot  \bar{X}_{j, d_{3}}^{0,\q} \rightarrow  \bar{X}_{j, d_{2}}^{0,-\q} \Leftrightarrow  \bar{X}_{i, d_{1}}^{0,\q} \cdot  X_{j, d_{3}}^{0,\q} \rightarrow  X_{j, d_{2}}^{0,-\q} $$

\begin{bfseries} Case 2: $d_{1} \geq 2^{2i-3}$  \end{bfseries}

\begin{bfseries} Subcase 2.1: $i \leq j-4-\log_{2}{n}$ \end{bfseries}

 Since there is a $\min$ involved on the right hand side we prove that we get a result with each of the terms there.

\begin{bfseries} Estimate for $2^{-i}$ \end{bfseries}

We derive this estimate by duality from {\mathversion{bold} $ X_{j, d_{2}}^{0,\q} \cdot  \bar{X}_{j, d_{3}}^{0,\q} \rightarrow  \bar{X}_{i, d_{1}}^{0,-\q}$}. Making use of (\ref{ge2}) for $f\delta_{P_{c_{1}}} * g \delta_{\bar{P}_{c_{2}}}$ and summing up on paraboloids, we obtain:

$$||B(u,v)||_{L^{2}} \leq 2^{\frac{(n-1)i}{2}} 2^{-\frac{j}{2}} ||u||_{X_{j, d_{2}}^{0,\q}} ||v||_{\bar{X}_{j, d_{3}}^{0,\q}}$$ 

Since $d_{1} \geq 2^{2i-3}$ we can conclude:

$$||B(u,v)||_{\bar{X}^{0,-\q}_{i,d_{1}}} \leq 2^{\frac{(n-3)i}{2}} 2^{-\frac{j}{2}} ||u||_{X_{j, d_{2}}^{0,\q}} ||v||_{\bar{X}_{j, d_{3}}^{0,\q}}$$ 

\begin{bfseries} Estimate for $d_3^{-\q}$ \end{bfseries}

We split $v=\sum_{l=1}^{n}v_l$ such that $\hat{v}_{l}$ is localized in a region where $|\xi_{l}|=\max_{1 \leq p \leq n} (|\xi_{p}|)$. For each $l$ we call $n_{l}$ the direction giving the coordinate $\xi_{l}$; for our purpose $n_l$ should be thought as a normal to the paraboloid in the region where $\hat{v}_{l}$ is localized. The main property we use is that:

$$||\hat{v}_{l}||_{L^{2}_{*}L^{1}_{n_{l}}} \les (2^{-j}d_2)^{\q} ||\hat{v}_{l}||_{L^{2}} \les 2^{-\frac{j}{2}}||v||_{X^{0,\q}_{j,d_{2}}}$$

\noindent
where $L^{2}_{*}$ is meant to be taken in the other $n-1 + 1$ directions (last one is $\tau$). Then we have:

$$||\hat{u} * \hat{v}_l||_{L^{2}} \leq ||\hat{u}||_{L^{1}_{*}L^{2}_{n_l}} ||\hat{v}_l||_{L^{2}_{*}L^{1}_{n_l}} \leq 2^{\frac{n+1}{2}i} 2^{-\frac{j}{2}}  ||\hat{u}||_{L^{2}} ||v_l||_{X^{0,\q}_{j,d_{2}}} \leq 2^{\frac{n-1}{2}i} 2^{-\frac{j}{2}} ||u||_{X^{0,\q}_{i,d_{1}}} ||v_l||_{X^{0,\q}_{j,d_{2}}}$$

The support of $\hat{v}_l$ (in the $\xi$ variable) is a solid section in the cone $|\xi_{l}|=\max_{1 \leq p \leq n} (|\xi_{p}|)$. Its sizes are $\approx 2^j$ in the direction giving the $\xi_l$ coordinate and $\approx \frac{2^j}{n} \geq 2^{i+4}$ in the others. Therefore the supports of $\hat{u} * \hat{v}_l$ are disjoint with respect to $l$. Summing with respect to $l$ and passing to $X^{0,\q}_{j,d_{3}}$, gives us the claimed estimate.

\begin{bfseries} Estimate for $d_2^{-\q}$ \end{bfseries}

We apply the same argument in Case 2 to {\mathversion{bold} $ \bar{X}_{i, d_{1}}^{0,\q} \cdot  X_{j, d_{2}}^{0,\q} \rightarrow  X_{j, d_{3}}^{0,-\q}$} and then by conjugation and duality we obtain the result we need.

The argument can be easily adapted for the estimates $\tilde{B}(\bar{u},v)$ and $\tilde{B}(u,\bar{v})$.

\begin{bfseries} Subcase 2.2: $i \geq j-4-\log_{2}{n}$ \end{bfseries}

 We obtain this estimate by duality from:

$$\bar{X}_{j, d_{3}}^{0,\q} \cdot  X_{j, d_{2}}^{0,\q} \rightarrow  \bar{X}_{i, d_{1}}^{0,-\q} \Leftrightarrow  X_{j, d_{3}}^{0,\q} \cdot  \bar{X}_{j, d_{2}}^{0,\q} \rightarrow  X_{i, d_{1}}^{0,-\q}$$

The last estimate is treated at the end of the proof, see the high-high to low interactions.

\vspace{.2 in}

{\mathversion{bold} $ X_{i, d_{1}}^{0,\q} \cdot   \bar{X}_{j, d_{2}}^{0,\q} \rightarrow  X_{j, d_{3}}^{0,-\q}$}

\vspace{.1 in}

If $j-4-\log_{2}{n} \leq i \leq j$ then the estimate is similar to the one in {\mathversion{bold} $ \bar{X}_{i, d_{1}}^{0,\q} \cdot  X_{j, d_{2}}^{0,\q} \rightarrow  X_{j, d_{3}}^{0,-\q}$}. If $i \leq j-4-\log_{2}{n}$, then we have the following cases:

\begin{bfseries} Case 1: $d_{2},d_{3} \leq 2^{2j-4}$ \end{bfseries}

This is incompatible since functions in $\bar{X}_{i, d_{1}}^{0,\q}$ have their Fourier transform supported in a region with $\tau < 0$ and functions in $X_{j, d_{3}}^{0,-\q}$ have their Fourier transform supported in a region with $\tau > 0$; an easy computation shows that, by convolution, the Fourier transform of a function in $X_{i, d_{1}}^{0,\q}$ cannot move the first support to the second one. 

\begin{bfseries} Case 2: $\max{(d_{2},d_{3})} \geq 2^{2j-3}$ \end{bfseries}

\begin{bfseries} Subcase 2.1: $d_{3} \geq 2^{2j-3} \geq d_{2}$ \end{bfseries}

This case can be treated the same way like {\mathversion{bold} $ X_{i, d_{1}}^{0,\q} \cdot  X_{j, d_{2}}^{0,\q} \rightarrow  X_{j, d_{3}}^{0,-\q}$}, just that we use this time the estimate for $f \delta_{P_{c_{1}}} * g \delta_{\bar{P}_{c_{2}}}$. Notice that the condition $d_{3} \geq 2^{2j-3}$ implies that we have to deal only with the estimate (\ref{b1}), since it becomes stronger than the estimate (\ref{b2}).

\begin{bfseries} Subcase 2.2: $d_{2} \geq 2^{2j-3} \geq d_{3}$ \end{bfseries}

This case can be obtained by duality from  {\mathversion{bold} $ X_{i, d_{1}}^{0,\q} \cdot  \bar{X}_{j, d_{3}}^{0,\q} \rightarrow  X_{j, d_{2}}^{0,-\q}$}, which is similar to the above estimate.

\begin{bfseries} Subcase 2.3: $\max{(d_{2},d_3)} \geq 2^{2j-3}$ \end{bfseries}

Use the $L^{1} * L^{2} \ra L^{2}$ argument. 

\vspace{.1 in}

\begin{bfseries} High - High interactions with output at low frequencies  \end{bfseries} 

\vspace{.1 in}

In order to be more suggestive about the kind of estimates we provide in this part of the proof, we choose to list them only for $i=j$. One would easily notice that we can replace one of the spaces on the left hand side with $X_{i, d_{1}}^{0,-\q}$ for any $|i-j| \leq 2$.

{\mathversion{bold} $  X_{i, d_{1}}^{0,\q} \cdot  X_{j, d_{2}}^{0,\q} \rightarrow   X_{k, d_{3}}^{0,-\q}$}

\vspace{.1 in}

Conjugation and duality give us:

$$ X_{k, d_{3}}^{0,\q} \cdot \bar{X}_{j, d_{2}}^{0,\q} \rightarrow  X_{j, d_{1}}^{0,-\q} \Rightarrow  \bar{X}_{k, d_{3}}^{0,\q} \cdot  X_{j, d_{2}}^{0,\q} \rightarrow  \bar{X}_{j, d_{1}}^{0,-\q} \Rightarrow X_{j, d_{1}}^{0,\q} \cdot X_{j, d_{2}}^{0,\q} \rightarrow  X_{k, d_{3}}^{0,-\q}$$

\noindent
and this is enough to justify the estimate. 

Nevertheless we owe a proof for the case: $k \geq j-4-\log_{2}{n}$ and $d_{3} \geq 2^{2k-3}$, see Subcase 2.2 from {\mathversion{bold} $  X_{i, d_{1}}^{0,\q} \cdot  X_{j, d_{2}}^{0,\q} \rightarrow   X_{j, d_{3}}^{0,-\q}$}.

\begin{bfseries} Case 1: $d_{1}, d_{2} \leq 2^{2j-2}$ \end{bfseries}

The argument is similar to Subcases 1.1 in the previous estimates. Making use of (\ref{ge}) we get:

$$||\hat{u} * \hat{v}||_{L^{2}} \les \int_{I_{1}} \int_{I_{2}} ||\hat{u} \delta_{P_{b_{1}2^{i}}} * \hat{v} \delta_{P_{b_{2}2^{j}}} ||_{L^{2}} 2^{-2j} db_{1} db_{2}$$

$$ \les \int_{I_{1}} \int_{I_{2}} 2^{(\frac{n}{2}-2)j} ||\hat{u}||_{L^{2}(P_{b_{1}2^{j}})} ||\hat{v}||_{L^{2}(P_{b_{2}2^{j}})} db_{1} db_{2} \les 2^{(\frac{n}{2}-1)j} ||u||_{X^{0,\q}_{j,d_{1}}} ||v||_{X^{0,\q}_{j,d_{2}}}$$

Next

$$||\tilde{B}(u,v)||_{X^{0,-\q}_{k,d_{3}}} \approx  (2^{2j})^{-\q} ||\tilde{B}(u,v)||_{L^{2}(A_{k,d_{3}})} \les 2^{(\frac{n}{2}-2)j} ||u||_{X^{0,\q}_{j,d_{1}}} ||v||_{X^{0,\q}_{j,d_{2}}}$$

\noindent
where we use the fact that $d_{3} \geq 2^{2k-3} \geq 2^{2j-7-2\log_{2}{n}}$. 

\begin{bfseries} Case 2: $\max{(d_{1}, d_{2})} \geq 2^{j-2}$ \end{bfseries}

This case is similar to Subcases 1.4 in the previous estimates and uses the trivial $L^{1} * L^{2} \rightarrow L^{2}$ argument. We skip the rest of the details. 

\vspace{.1 in}

{\mathversion{bold} $X_{j, d_{1}}^{0,\q} \cdot  \bar{X}_{j, d_{2}}^{0,\q} \rightarrow   X_{k, d_{3}}^{0,-\q}$}

\vspace{.1 in}

In the same way as above, duality gives us the estimates as claimed in the Theorem. 

\end{proof}
 
We are done with all the preparatory results for this section and we turn our attention to the proof of the main result. 
 
\begin{proof}[Proof of Theorem \ref{bil}]

The proof is a direct consequence of the estimates in Proposition \ref{u4} since we only have to perform the summations with respect to the modulations $d$. For the low-high to high interactions we can simplify the model and consider $k=j$. Also working with $u \in X^{0,\q,1}_i$ is, essentially, the same with working with $u_i \in X^{0,\q,1}$. We prefer the second one, due to the second index coming from the modulation, namely $u_{i,d}$. 

We decompose:

\beq \label{dec10}
B(u_{i},v_{j})_j = \sum_{d_{1},d_2,d_3} B(u_{i,d_1},v_{j,d_2})_{j,d_3}=
\eeq

$$\sum_{\max{(d_2,d_3)} \leq 2^{2i}} \sum_{d_1}B(u_{i,d_1},v_{j,d_2})_{j,d_3} + \sum_{ 2^{2i} \leq  \max{(d_2,d_3)} \leq 2^{i+j+5}}\sum_{d_1} B(u_{i,d_1},v_{j,d_2})_{j,d_3}+$$

$$\sum_{ 2^{i+j+6} \leq  \max{(d_2,d_3)}} \sum_{d_1} B(u_{i,d_1},v_{j,d_2})_{j,d_3}$$

We take separately each term and estimate it by using the results in Proposition \ref{u4}. For the first term ($\max{(d_2,d_3)} \leq 2^{2i}$), we use \eqref{b100}:

$$||B(u,v)||_{X^{0,-\q}_{j,d_{3}}} \les 2^{\frac{n+1}{2}i} 2^{-\frac{j}{2}} ||u||_{X^{0,\q}_{i,d_{1}}} ||v||_{X^{0,\q}_{j,d_{2}}}$$

Summing this estimates with respect to $d_1$ and $d_2,d_3$ (up to $2^{2i}$) we obtain:

$$||B(u_{i},v_{j, \leq 2^{2i}})||_{X^{0,-\q,1}_{j, \leq 2^{2i}}} \les i^2 2^{\frac{n}{2}i} ||u_i||_{X^{0,\q}_{i}} ||v_{j, \leq 2^{2i}}||_{X^{0,\q}}$$

For the second term in (\ref{dec10}), we notice that there is a symmetry in $d_2,d_3$ in the estimates coming from Proposition \ref{u4}. Let us assume that $d_{2} \leq d_{3}$. Then \eqref{b100} gives us the estimate:

$$||B(u,v)||_{X^{0,-\q}_{j,d_{3}}} \les 2^{\frac{n-1}{2}i} 2^{-\frac{j}{2}} d_{3}^{\q} ||u||_{X^{0,\q}_{i,d_{1}}} ||v||_{X^{0,\q}_{j,d_{2}}}$$

In the case $d_2 \geq d_3$ we obtain:

$$||B(u,v)||_{X^{0,-\q}_{j,d_{3}}} \les 2^{\frac{n-1}{2}i} 2^{-\frac{j}{2}} d_{2}^{\q} ||u||_{X^{0,\q}_{i,d_{1}}} ||v||_{X^{0,\q}_{j,d_{2}}}$$

Notice that in both cases we can bound the coefficient by $2^{\frac{n}{2}i}$, due to the condition $\max{(d_2,d_3)} \leq 2^{i+j+5}$. Summing up with respect to $d_1,d_2,d_3$ in the corresponding range we obtain:

$$||\sum_{ 2^{2i} \leq  \max{(d_2,d_3)} \leq 2^{i+j+5}}\sum_{d_1} B(u_{i,d_1},v_{j,d_2})_{j,d_3}||_{X^{0,-\q,1}} \les j^2 2^{\frac{n}{2}i} ||u||_{X^{0,\q,1}_{i}} ||v||_{X^{0,\q,1}_{j}}$$

We should remark that even a more careful use of the above estimates would not completely eliminate the logarithmic of the high frequency. At the end of the argument we will slightly modify the range of summation such that we obtain a more convenient estimate. 

For the case $\max{(d_2,d_3)} \geq 2^{i+j+6}$ we use (\ref{b102}) together with the observation there and claim:

$$||B(u,v)||_{X^{0,-\q,1}_{j}} \les  2^{\frac{n}{2}i} ||u||_{X^{0,\q,1}_{i}} ||v||_{X^{0,\q,1}_{j, \geq 2^{i+j+6}}}$$

Bringing together all the estimates we obtain for the sums in the decomposition (\ref{dec10}) and then passing to general $s$ gives us:

$$||B(u,v)||_{X^{s,-\q,1}_{j}} \les  j^2 2^{(\frac{n}{2}-s)i} ||u||_{X^{s,\q,1}_{i}} ||v||_{X^{s,\q,1}_{j}}$$

Going back to the estimates we observe that if we want to cover only the range of $d_2,d_3 \geq 2^{j-i}$, then we pick up only logarithms of the low frequency:

$$||B(u,v)||_{X^{s,-\q,1}_{j, \geq 2^{j-i}}} \les  i^2 2^{(\frac{n}{2}-s)i} ||u||_{X^{s,\q,1}_{i}} ||v||_{X^{s,\q,1}_{j,\geq 2^{j-i}}}$$

Based on (\ref{b103}) we can run the same argument for estimating the high-high to low frequencies and obtain:

$$||B(u,v)||_{X^{s,-\q,1}_{k}} \les  k j 2^{(\frac{n}{2}-1+s)(k-j)} 2^{(\frac{n}{2}-s)j} ||u||_{X^{s,\q,1}_{i}} ||v||_{X^{s,\q,1}_{j}} $$

\noindent
where $|i-j| \leq 2$.

\end{proof}

\section{Bilinear estimates involving $Y$ spaces} \label{bes+}

This section completes the theory of bilinear estimates. As discussed in the previous section, the estimates provided by Theorem \ref{bil} are not satisfactory for all dyadic ranges. We decided that if we analyze the interactions of type $B(u_i,v_j)$ with $i \leq j$, then we use the $X^{s,\q,1}$ spaces as long as $2(n+2)i \geq j$. Otherwise we have to complete the theory of bilinear estimates by involving the the more delicate structure $Y$.

Throughout this section we assume that $2(n+2)i \leq j$. The main result is the following theorem. 

\begin{t1} \label{yt} We have the estimate

\beq \label{y3}
||B(u_{i}, v_{j, \leq 2^{i+j+5}})||_{W_{j,\leq 2^{i+j+5}}} \leq i^{2} 2^{\frac{n}{2}i} ||u_{i}||_{Z} ||v_{j}||_{Z}
\eeq

 As a direct consequence we obtain:

\beq \label{y6}
||B(u_{i}, v_{j})||_{W_{j}^s} \leq i^{2} 2^{(\frac{n}{2}-s)i} ||u_{i}||_{Z^s} ||v_{j}||_{Z^s}
\eeq

\end{t1}

\vspace{.1in}

We need some preparatory results. In what follows one should think of $g$ as being either $g_{\xi,\xi^{2}+k}$ or $g_{\xi,l}$ (recall the concepts related to the decompositions in \eqref{aux2} and \eqref{aux1}).

We assume $\hat{g}$ is localized in frequency on a scale $2^{-i} \times ... \times 2^{-i} \times 1$ ($\xi \times \tau$), hence the dual scale to localize in the physical space is $2^{i} \times ... \times 2^{i} \times 1$ ($x \times t$). We define the system of cubes $(Q^{m}_{i})_{m \in \Z^{n}}$ to be a partition of $\R^{n}$ with the properties: $Q^{m}_{i}$ is centered at $2^{i}m$ and has sizes $2^{i} \times ... \times 2^{i}$. Associated to this system we build a partition $(Q_{i}^{m,l})_{(m,l) \in Z^{n+1}}$ of $\R^{n+1}$ defined by:
$$Q_{i}^{m,l}=\cup_{t \in [l,l+1]} Q^{m}_{i} \times \{ t \}= Q^{m}_{i} \times [l,l+1]$$

We have the following refinement of the Sobolev embedding:   
 
\begin{l1}

Let $g \in L^{2}$ such that $\hat{g}$ is supported in a tube of size $2^{-i} \times ... \times 2^{-i} \times 1$. We have the estimate:

\beq \label{c0}
\sum_{(m,l) \in Z^{n+1}} ||g||^{2}_{L^{\infty}(Q^{m,l}_{i})} \les 2^{-ni} ||g||^{2}_{L^{2}}
\eeq

\end{l1}

This result says that we obtained the standard $L^{\infty}$ estimates on tubes localized in the physical side on the dual scale and then sum the estimates in $l^2$ with respect to the tubes.  

Using this result we obtain 

\begin{p1} Let $f$ and $g$ be two functions with the following properties: $f=f_{\eta} \in Y_{j}$, $|\eta| \approx 2^{j}$, $g \in \mathcal{D}L^{2}$, $\hat{g}$ is supported at frequency $2^{i}$ in a tube of size $2^{-i} \times ... \times 2^{-i} \times 1 (\xi \times \tau)$. Then we have the estimates:

\beq \label{m9}
|| f \cdot  g||_{L^{2}} \les  2^{-\frac{(n-1)i+j}{2}} ||f||_{Y_{j}} ||g||_{L^{2}}
\eeq

\end{p1}

\begin{proof}

 For each $(m,l) \in \Z^{3}$ we denote by $C^{m,l}=\{m' \in Z^{2}: Q^{m',l} \cap T^{m,l}_{\eta} \ne \emptyset \}$. The intersection $Q^{m',l} \cap T^{m,l}_{\eta}$, if nonempty, has size at most $2^{i-j}$ in the direction of $t$. We have:

$$||f \cdot g||^{2}_{L^{2}(T^{m,l}_{\eta})}=\sum_{m' \in C^{m,l}} ||f \cdot g||^{2}_{L^{2}(Q_{i}^{m',l} \cap T^{m,l}_{\eta})} \les$$

$$ 2^{i-j} \sum_{m' \in C^{m,l}} ||f||^{2}_{L_{t}^{\infty}L^{2}_{x}(T^{m,l}_{\eta})} ||g||^{2}_{L^{\infty}(Q_{i}^{m',l})} \les $$

$$2^{i-j} ||f||^{2}_{L_{t}^{\infty}L^{2}_{x}(T^{m,l}_{\eta})} \sum_{m' \in C^{m,l}} ||g||^{2}_{L^{\infty}(Q_{i}^{m',l})} \les   2^{-(n-1)i-j} ||f||^{2}_{L_{t}^{\infty}L^{2}_{x}(T^{m,l}_{\eta})} ||g||^{2}_{L^{2}}$$

In the last estimate we have used the result in (\ref{c0}).

We sum the above estimate with respect to $(m,l)$ over $Z^{3}$ to obtain (\ref{m9}).

\end{proof}

On behalf of (\ref{m9}) we derive estimates adapted to our spaces.

\begin{p1}

 If $d \leq 2^{j-i}$ we have the estimates:

\beq \label{y1}
||\tilde{B}(u_{i},v_{j, \leq d})||_{L^{2}} \les 2^{\frac{(n-1)i-j}{2}} ||u_{i}||_{X^{0,\q,1}} ||v_{j, \leq d}||_{Y^0}
\eeq

\beq \label{y2}
||\tilde{B}(u_{i},v)||_{\mathcal{Y}^0_{j, \leq d}} \les 2^{\frac{(n-1)i-j}{2}} ||u_{i}||_{X^{0,\q,1}} ||v||_{L^{2}}
\eeq

The above estimates hold true if $u_i$ is replaced by $\bar{u}_i$.

\end{p1}

\begin{r2} The condition $2(n+2)i \leq j$ implies that the high frequency does not see the difference between $X_i^{0,\q,1}$ and $\bar{X}_i^{0,\q,1}$. Therefore once we obtain the estimates for $u_{i}$, we can easily deduce the ones for $\bar{u}_i$.
\end{r2}

\begin{proof}

(\ref{y1}) and (\ref{y2}) are dual to each other. We choose to prove (\ref{y1}) in the particular case $d=2^{j-i}$. The proof can be easily adapted for all other $d$'s.

We decompose $\R^{n+1}$ is disjoint parallelepipeds of sizes $2^{i} \times ... \times 2^i \times 2^{2i}$ ($n$ $\xi$ directions $\times \tau$). We intersect the set $A_{j,\leq 2^{j-i}} = \{(\xi,\tau)|: |(\xi,\tau)| \approx 2^{i}, |\tau-\xi^{2}| \leq 2^{j-i}\}$ with these parallelepipeds to obtain a decomposition of $A_{j,\leq 2^{j-i}}$. Note that for each parallelepiped $R$ intersecting $A_{j,\leq 2^{j-i}}$ there is a direction in which $R \cap A_{j,\leq 2^{j-i}}$ has size $2^{j-i}$ - one could think of this direction as the projection of the normal to $P$ at one of the points in the intersection onto the $\xi$ space. We name this vector by $n_{R}$. This way we obtain the decomposition

 $$v_{j,\leq 2^{-i}}=\sum_{R} v_{R}$$

\noindent 
where $\hat{v}_R$ is the part of $\hat{v}_{j,\leq 2^{-i}}$ supported in $R \cap A_{j,\leq 2^{j-i}}$. The support of $\hat{u}_i$ is included in $R_0$ (the parallelepiped containing the origin), therefore the supports of $\hat{u}_i * \hat{v}_{R}$ are disjoint with respect to $R$. As a consequence

$$||u_{i} \cdot v_{j,\leq 2^{-i}}||_{L^{2}}^{2} \approx \sum_{R} ||u_{i} \cdot v_{R}||_{L^{2}}^{2}$$

Hence we should estimate $||u_{i} \cdot v_{R}||_{L^{2}}$ for each $R$. We continue with an obvious estimate:

$$||u_{i} \cdot v_{R}||_{L^{2}} \les \sum_{d} ||u_{i,d} \cdot v_{R}||_{L^{2}}$$

Next we distinguish two cases according to whether $d \leq 2^{2i-4}$ or $d \geq 2^{2i-3}$. 

Let us treat the first one, $d \leq 2^{2i-4}$. We decompose:

$$u_{i,d} = \sum_{k=2^{-1}d}^{2d} \sum_{\xi \in \Xi^{i}} u_{\xi,\xi^{2}+k}$$

 The key property is that for fixed $k$ the support $\hat{u}_{\xi,\xi^{2}+k} * \hat{v}_{R}$ is disjoint with respect to $\xi$, as we vary $\xi$ in the direction of $n_{R}$. One can do an explicit computation in this direction; we prefer instead to provide a more intuitive one. As we move $\xi$ in the direction of $n_{R}$, we move the support of $\hat{u}_{\xi,\xi^{2}+k} * \hat{v}_{R}$ in the normal direction to $P$ with at least $2^{-i}$. We also move the support in the direction of $\tau$ with at most $2^{2i}$ which accounts for a correction of the distance to $P$ of $\approx 2^{2i-j} \leq 2^{-2i}$, hence insignificant to the original shift. Since the support of  $\hat{u}_{\xi,\xi^{2}+k} * \hat{v}_{R}$ is included in a strip of width $2^{-i}$ in the direction of $n_{R}$ our conclusion follows. 

We lose orthogonality of the interaction  $\hat{u}_{\xi,\xi^{2}+k} * \hat{v}_{R}$ as we move $\xi$ in the transversal directions to $n_{R}$; from the total of $2^{2ni}$ values for $\xi \in \Xi^{i}$, $2^{2i}$ go in the direction of $n_{R}$ and $2^{2(n-1)i}$ in the other ones. Hence we can conclude:

$$||\sum_{\xi} u_{\xi,\xi^{2}+k} \cdot v_{R} ||_{L^{2}} \les 2^{(n-1)i} ( \sum_{\xi} || u_{\xi,\xi^{2}+k} \cdot v_{R} ||^{2}_{L^{2}} )^{\q}$$

For fixed $\xi$ we have:

$$|| u_{\xi,\xi^{2}+k} \cdot v_{R} ||^{2}_{L^{2}} \approx \sum_{\eta} || u_{\xi,\xi^{2}+k} \cdot v_{\eta, \leq 2^{j-i}} ||^{2}_{L^{2}}$$

The second summation is performed over the range of those $\eta$ for each the support of $\hat{v}_{\eta, \leq 2^{j-i}}$ intersects $R$. We have reduced matters to the use of the estimate in (\ref{m9}), hence we obtain:

$$||u_{\xi,\xi^{2}+k} \cdot v_{\eta, \leq 2^{j-i}} ||_{L^{2}} \les 2^{-\frac{(n-1)i+j}{2}} ||u_{\xi,\xi^{2}+k}||_{L^2}  || v_{\eta, \leq 2^{j-i}} ||_{Y^0}$$

Going reverse in the above argument we obtain:

$$||\sum_{\xi} u_{\xi,\xi^{2}+k} \cdot v_{R} ||_{L^{2}} \les 2^{\frac{(n-1)i-j}{2}} ||\sum_{\xi} u_{\xi,\xi^{2}+k}||_{L^2}  || v_{R} ||_{Y^0}$$

Then we sum up with respect to $k \in [2^{-1}d, 2d]$ gives us:

$$||u_{i,d} \cdot v_{R} ||_{L^{2}} \les 2^{\frac{(n-1)i-j}{2}} ||u_{i,d}||_{X^{0,\q}}  || v_{R} ||_{Y^0}$$

In order to continue with the full argument we have to deal with the case $d \geq 2^{2i-3}$.  In principle the approach is similar; we decompose

$$u_{i,d} = \sum_{k} \sum_{\xi \in \Xi^{i}} u_{\xi,k}$$

\noindent
and then fix $k$ and rewrite the same argument as before to claim:

$$||u_{i,d} \cdot v_{R} ||_{L^{2}} \les 2^{\frac{(n-1)i-j}{2}} ||u_{i,d}||_{X^{0,\q}}  || v_{R} ||_{Y^0}$$

Now we can perform the summation with respect to $d$ and then the one with respect to $R$'s to claim

$$||u_{i} \cdot v_{j,\leq 2^{j-i}}||_{L^2} \les 2^{\frac{(n-1)i-j}{2}} ||u_{i}||_{X^{0,\q,1}} ||v_{j,\leq 2^{j-i}}||_{Y^0}$$

\end{proof}

\begin{proof}[Proof of Theorem \ref{yt}] 

 We tacitly agree that in the sums bellow we have the bounds $d_2,d_3 \leq 2^{i+j+5}$ and decompose

\beq \label{s1}
\sum_{d_2,d_3} B(u_i,v_{j,d_2})_{j,d_3} = \sum_{d_2 \leq d_3} B(u_i,v_{j,d_2})_{j,d_3} +  \sum_{d_2 > d_3} B(u_i,v_{j,d_2})_{j,d_3}=
\eeq

$$\sum_{d_3} B(u_i,v_{j,\leq d_3})_{j,d_3} + \sum_{d_2} B(u_i,v_{j,d_2})_{j, < d_{2}}$$

We decompose the first term in the sum:

$$\sum_{d_3} B(u_i,v_{j,\leq d_3})_{j,d_3}= \sum_{d_3 \leq 2^{j-i}} B(u_i,v_{j,\leq d_3})_{j,d_3} + \sum_{d_3 \geq 2^{j-i}} B(u_i,v_{j,\leq d_3})_{j,d_3}=$$

$$\sum_{d_3 \leq 2^{j-i}} B(u_i,v_{j,\leq d_3})_{j,d_3} + \sum_{d_3 \geq 2^{j-i}} \left( B(u_i,v_{j,\leq 2^{j-i}})_{j,d_3} + \sum _{d=2^{j-i}}^{d_{3}} B(u_i,v_{j, d})_{j,d_3}\right)$$

Using (\ref{r1}), (\ref{res1}), (\ref{y1}) and the fact that $2(n+2)i \leq j$, we estimate the first sum as follows:

$$||B(u_i,v_{j,\leq d_3})||_{X^{0,-\q}_{j,d_{3}}} \les d^{-\q}_{3} \max{(2^{2i},d_3)} 2^{\frac{(n-1)i-j}{2}} ||u_{i}||_{X^{0,\q,1}} ||v_{j,\leq d_3}||_{Y^0} \les$$

$$ 2^{\frac{n}{2}i} \max{((2^{-i-j}d_{3})^{\q}, 2^{-\frac{j}{10}})} ||u_{i}||_{X^{0,\q,1}} ||v_{j,\leq d_3}||_{Y^0}$$

For the second term in the sum, we obtain in a similar way:

$$||B(u_i,v_{j,\leq 2^{j-i}})||_{X^{0,-\q}_{j,d_{3}}} \leq  2^{i} \max{((2^{-i-j}d_{3})^{\q}, 2^{-\frac{j}{10}})}||u_{i}||_{X^{0,\q,1}} ||v_{j,\leq 2^{j-i}}||_{Y^0}$$

Taking into account that:

$$\sum_{d_3 \leq 2^{i+j+5}} \max{((2^{-i-j}d_{3})^{\q}, 2^{-\frac{j}{10}})} \les 1$$

\noindent
and that the terms of type $ ||B(u_i,v_{j, d})||_{X_{j,d_3}^{0,-\q}}$ for $d,d_{3} \geq 2^{j-i}$ have been chosen to be treated in $X^{0,\q,1}$ we conclude with:

$$||\sum_{d_3} B(u_i,v_{j,\leq d_3})||_{X^{0,-\q,1}_{j,d_3}} \leq i^{2} 2^{\frac{n}{2}i} ||u_{i}||_{Z^0} ||v_{j}||_{Z^0}$$

For the second sum in (\ref{s1}) we proceed in a similar manner:

$$\sum_{d_2} B(u_i,v_{j, d_2})_{j, < d_2}= \sum_{d_2 \leq 2^{j-i}} B(u_i,v_{j, d_2})_{j, < d_2} + \sum_{d_2 \geq 2^{j-i}} B(u_i,v_{j, d_2})_{j, < d_2}=$$

$$\sum_{d_2 \leq 2^{j-i}} B(u_i,v_{j, d_2})_{j, < d_2} + \sum_{d_2 \geq 2^{j-i}} \left( B(u_i,v_{j, d_2})_{j, \leq 2^{j-i} } + \sum _{d=2^{j-i}}^{d_{2}} B(u_i,v_{j, d_2})_{j,d}\right)$$

Then we continue the same way as we did with the first term in the sum, just that we use this time (\ref{y2}). We obtain:

$$||\sum_{d_2 \leq 2^{j-i}} B(u_i,v_{j, d_2})_{j, < d_2} + \sum_{d_2 \geq 2^{j-i}} B(u_i,v_{j, d_2})_{j, \leq 2^{j-i} }||_{\mathcal{Y}^0} \leq 2^{\frac{n}{2}i} ||u_{i}||_{Z^0} ||v_{j}||_{Z^0}$$

The terms $B(u_i,v_{j, d_2})_{j,d}$, for $d,d_2 \geq 2^{j-i}$ have been chosen from before to be treated in $X^{0,\q,1}$ (see (\ref{be22})), hence we can conclude with:

$$||\sum_{d_2} B(u_i,v_{j, d_2})_{j, < d_2}||_{W^0} \leq i^{2} 2^{\frac{n}{2}i} ||u_{i}||_{Z^0} ||v_{j}||_{Z^0}$$

This ends the argument for (\ref{y3}). 

\eqref{y6} is a direct consequence of the estimate in part \eqref{y3} and of the estimate in (\ref{b2}).

\end{proof}

We finish this section with the proof of the key estimate for our problem, namely the bilinear estimate.

\begin{proof}[Proof of (\ref{be})]

We decompose:

$$B(u,v)=\sum_{i,j,k} B(u_i,v_j)_k$$

Taking into account only the nontrivial interactions we can write:

\beq \label{dec15}
B(u,v)=\sum_{|i-j| \geq 2} \sum_{|k-\max{(i,j)}| \leq 2} B(u_i,v_j)_k + \sum_{|i-j| \leq 1} \sum_{k \leq \max{(i,j)+2}} B(u_i,v_j)_k
\eeq

Due to the symmetry of the indexes, it is enough to estimate the first sum above for the particular case $i \leq j-2$. Using (\ref{be11}) and (\ref{y6}) we obtain:

$$|| \sum_{|j-k| \leq 2} \sum_{i \leq j-2} B(u_i,v_j)_k||^2_{W^s} \les $$

$$ \sum_{|j-k| \leq 2} ||\sum_{\frac{j}{2(n+2)} \leq i \leq j-2} B(u_i,v_j)_k||^{2}_{X^{s,-\q,1}} + \sum_{|j-k| \leq 2} ||\sum_{i \leq \frac{j}{2(n+2)}} B(u_i,v_j)_k||^{2}_{W^s} \les$$

$$ \sum_{|j-k| \leq 2} \left( \sum_{\frac{j}{2(n+2)} \leq i \leq j-2}  k^2 2^{(\frac{n}{2}-s)i} ||u_i||_{X^{s,\q,1}} ||v_{j}||_{X^{s,\q,1}} \right)^2 + $$

$$\sum_{|j-k| \leq 2} \left( \sum_{i \leq \frac{j}{2(n+2)}} i^2 2^{(\frac{n}{2}-s)i} ||u_i||_{X^{s,\q,1}} ||v_{j}||_{Z^{s}} \right)^2 \les $$

$$ \sum_{|j-k| \leq 2} \left( \sum_{\frac{j}{2(n+2)} \leq i \leq j-2}  k^2 i j 2^{(\frac{n}{2}-s)i} ||u_i||_{Z^{s}} ||v_{j}||_{Z^{s}} \right)^2 +$$

$$ \sum_{|j-k| \leq 2} \left( \sum_{i \leq \frac{j}{2(n+2)}} i^3 2^{(\frac{n}{2}-s)i} ||u_i||_{Z^{s}} ||v_{j}||_{Z^{s}} \right)^2 \les $$

$$C_s^2 \sum_{|j-k| \leq 2}  ||u||^2_{Z^{s}} ||v_{j}||^{2}_{Z^{s}} \les C_s^2 ||u||^2_{Z^{s}} ||v||^{2}_{Z^{s}}$$

We have used the fact that $\frac{n}{2} < s$. For the second term in the decomposition in (\ref{dec15}) we use (\ref{be11})

$$||\sum_{|i-j| \leq 1} \sum_{k \leq \max{(i,j)+2}} B(u_i,v_j)_k||_{W^s} \les ||\sum_{|i-j| \leq 1} \sum_{k \leq \max{(i,j)+2}} B(u_i,v_j)_k||_{X^{s,-\q,1}} \les$$

$$ \sum_{|i-j| \leq 1} \sum_{k \leq \max{(i,j)+2}} j^2 2^{(\frac{n}{2}-s)j} 2^{(\frac{n}{2}-1+s)(k-j)}||u_i||_{X^{s,\q,1}} ||v_j||_{X^{s,\q,1}} \les$$

$$ \sum_{|i-j| \leq 1}  j^2 2^{(\frac{n}{2}-s)j} ||u_i||_{X^{s,\q,1}} ||v_j||_{X^{s,\q,1}} \les \sum_{|i-j| \leq 1}  j^4 2^{(\frac{n}{2}-s)j} ||u_i||_{Z^{s}} ||v_j||_{Z^{s}} \les$$

$$C_s ||u||_{Z^s} ||v||_{Z^s}$$

In the last line we have used that $\frac{n}{2} < s$ which implies $j^3 2^{(\frac{n}{2}-s)j} \leq C_s$ followed by a Cauchy-Schwarz. 

Bringing together the two estimates we obtained for the two sums in the decomposition (\ref{dec15}) gives us the claim (\ref{be}). 

\end{proof}

\section{Algebra Properties} \label{alg}

In this section we intend to prove (\ref{A}) and (\ref{M}). One has to keep in mind that in some sense these are weaker estimates than the bilinear ones. First we want to derive the equivalent of  (\ref{A}) in $X^{s, \q,1}$ and then prove it for $Z^{s}$. We discover that for $s > \frac{n}{2}$

\beq \label{y4}
 X^{s,\q,1}+\bar{X}^{s,\q,1} \ \mbox{is an complex algebra}.
\eeq

 In the case of the wave equation it was known that for $s > 1$ the corresponding $X^{s,\q,1}$ is an algebra. In our case we have to adjust the structure to $X^{s,\q,1}+\bar{X}^{s,\q,1}$ due to complex nature of the equation.  
 
 Once we understand this result we are ready to adjust the argument and prove (\ref{A}), i.e. that $Z^s+\bar{Z}^s$ is a complex algebra for $s > \frac{n}{2}$. 

Then by duality, conjugation and additional arguments (when needed) we derive the property (\ref{M}).

The structure of this section goes as described above. We prove (\ref{y4}) and (\ref{A}), and then we conclude with the proof of (\ref{M}).

\begin{proof}[Proof of (\ref{y4})] For the beginning, we assume that $u,v  \in X^{s,\q,1}$ and intend to show that $u v  \in X^{s,\q,1}$.  The argument we provide here is very similar to the one in the proof of Theorem \ref{bil}. We decompose:

\beq \label{y14}
\tilde{B}(u,v)= \sum_{i,j,k} \tilde{B}(u_{i},v_{j})_k
\eeq

If $i \leq j$ and $|k-j| \leq 2$, we can reduce the matters to the case $k=j$ and then start with the same decomposition as in (\ref{dec10})

$$\tilde{B}(u_{i},v_{j})_j = \sum_{d_{1},d_2,d_3} \tilde{B}(u_{i,d_1},v_{j,d_2})_{j,d_3}$$

$$ \sum_{ \max{(d_2,d_3)} \leq 2^{i+j+5}}\sum_{d_1} \tilde{B}(u_{i,d_1},v_{j,d_2})_{j,d_3}+ \sum_{ 2^{i+j+6} \leq  \max{(d_2,d_3)}} \sum_{d_1} \tilde{B}(u_{i,d_1},v_{j,d_2})_{j,d_3}$$

We take separately each term and estimate it by using the results in Proposition \ref{u9}. We rewrite (\ref{b99}) as 

$$||\tilde{B}(u,v)||_{X^{0,\q}_{j,d_{3}}} \les 2^{\frac{n}{2}i} (2^{-i-j} d_{3})^{\q}  ||u||_{X^{0,\q}_{i,d_{1}}} ||v||_{X^{0,\q}_{j,d_{2}}}$$

Summing this estimates with respect to $d_1$(up to $2^{2i+2}$) and with $d_2,d_3$ (up to $2^{i+j+4}$) we obtain:

$$||\tilde{B}(u_{i},v_{j, \leq 2^{i+j+4}})||_{X^{0,\q,1}_{j, \leq 2^{i+j+5}}} \les  2^{\frac{n}{2}i} ||u||_{X^{0,\q,1}_{i}} ||v||_{X^{0,\q,1}_{j, \leq 2^{i+j+5}}}$$

For the second term in the above sum, we recall the observation from (\ref{b96}) that unless $4^{-1} \leq d_{2}d_{3}^{-1} \leq 4$ we have a trivial estimate. Then using (\ref{b96}) we can sum with respect to all $d$'s to obtain

$$||\sum_{ 2^{i+j+6} \leq  \max{(d_2,d_3)}} \sum_{d_1} \tilde{B}(u_{i,d_1},v_{j,d_2})_{j,d_3}||_{X^{0,\q,1}} \les 2^{\frac{n}{2}i} ||u||_{X^{0,\q,1}_{i}} ||v||_{X^{0,\q,1}_{j}}$$

We put the two estimates together, pass to general $s$ and conclude with:

$$||\tilde{B}(u_{i},v_{j})||_{X^{s,\q,1}_{j}} \les  2^{(\frac{n}{2}-s)i} ||u||_{X^{s,\q,1}_{i}} ||v||_{X^{s,\q,1}_{j}}$$

We obtain the same result for the case $j \leq i$ and $|k-i| \leq 2$. For the high-high to low interactions, i.e. $|i-j| \leq 2$, we proceed as follows

$$\tilde{B}(u_{i},v_{j})_k = \sum_{d_{1},d_2,d_3} \tilde{B}(u_{i,d_1},v_{j,d_2})_{k,d_3}$$

We use (\ref{b95}) to obtain

$$||\tilde{B}(u,v)||_{X^{0,\q,1}_{k}} \les \sum_{d_{1},d_2,d_3} ||\tilde{B}(u_{i,d_1},v_{j,d_2})_{k,d_3}||_{X^{0,\q}_{k}}$$

$$ \sum_{d_{1},d_2,d_3} 2^{(\frac{n}{2}-2)k} d_{3}||u||_{X^{0,\q}_{i,d_{1}}} ||v||_{X^{0,\q}_{j,d_{2}}} \les 2^{\frac{n}{2}k}||u||_{X^{0,\q,1}_{i}} ||v||_{X^{0,\q,1}_{j}}$$

In the last line we have used the fact that the range of $d_3$ goes up to $ \approx 2^{2k}$. For general $s$ this becomes:

\beq \label{y11}
||\tilde{B}(u,v)||_{X^{s,\q,1}_{k}} \les 2^{(\frac{n}{2}+s)k} 2^{-2js} ||u||_{X^{s,\q,1}_{i}} ||v||_{X^{s,\q,1}_{j}}
\eeq

As we did before, see the proof of (\ref{be}), we can add up all these estimates and, for $s > 1$, obtain:

$$||\tilde{B}(u,v)||_{X^{s,\q,1}} \leq C_{s} ||u||_{X^{s,\q,1}} ||v||_{X^{s,\q,1}}$$

This corresponds to the assertion

$$X^{s,\q,1} \cdot X^{s,\q,1} \ra X^{s,\q,1}$$

By conjugation we get for free that

$$\bar{X}^{s,\q,1} \cdot \bar{X}^{s,\q,1} \ra \bar{X}^{s,\q,1}$$

We are left now with proving 

$$\bar{X}^{s,\q,1} \cdot X^{s,\q,1} \ra X^{s,\q,1}+ \bar{X}^{s,\q,1}$$

 If $u,v \in X^{s,\q,1}$ we have to show that $\bar{u} \cdot v \in X^{s,\q,1}+ \bar{X}^{s,\q,1}$. We split

$$\bar{u} \cdot v=\sum_{i,j} \bar{u}_{i} \cdot v_{j} = \sum_{i \leq j} \bar{u}_{i} \cdot v_{j}+\sum_{i >j} \bar{u}_{i} \cdot v_{j}$$

The argument we provided above works the same if we conjugate the low frequency, since all the estimates in Proposition \ref{u9} allow us to put a conjugate on the low frequency. Therefore we can claim

$$\sum_{i \leq j} \bar{u}_{i} \cdot v_{j} \in X^{s,\q,1}$$

Now we get for free (by conjugation) that $\sum_{i > j} \bar{u}_{i} \cdot v_{j} \in \bar{X}^{s,\q,1}$ and this ends our proof.

\end{proof}

\begin{proof}[Proof of (\ref{A})] We start by showing that if $u,v \in Z^s$ then $u \cdot v \in Z^s$. We use the decomposition in (\ref{y14}).

If $i \leq j$ and $|k-j| \leq 2$, we can reduce the matters to the case $k=j$. We further split $v_j= v_{j, \leq 2^{j-i}}+ v_{j, \geq 2^{j-i}}$. Since 

$$||v_{j, \geq 2^{j-i}}||_{X^{s,\q,1}} \leq i ||v_{j, \geq 2^{j-i}}||_{Z^{s}}$$

\noindent
we can invoke the arguments above to claim

$$||\tilde{B}(u_i,v_{j, \geq 2^{j-i}})||_{X^{s,\q,1}_{j}} \les i^2 2^{(\frac{n}{2}-s)i} ||u_i||_{Z^s} ||v_{j, \geq 2^{j-i}}||_{Z^{s}}$$

We observe that $\mathcal{F}(\tilde{B}(u_i,v_{j, \leq 2^{j-i}}))$ is localized in a region with $|\tau - \xi^{2}| \leq 2^{i+j+5}$. Using (\ref{y1}) 
we continue with

$$||\tilde{B}(u_i,v_{j, \leq 2^{j-i}})||_{X^{s,\q,1}_{j}} \leq \sum_{d \leq 2^{i+j+5}} ||\tilde{B}(u_i,v_{j, \leq 2^{j-i}})||_{X^{s,\q,1}_{j,d}} \les $$

$$\sum_{d \leq 2^{i+j+5}} 2^{js} d^{\q} ||\tilde{B}(u_i,v_{j, \leq 2^{j-i}})||_{L^{2}} \les 2^{js} \sum_{d \leq 2^{i+j+5}} 2^{\frac{(n-1)i-j}{2}} d^{\q} ||u_i||_{X^{0,\q,1}} ||v_{j, \leq 2^{j-i}}||_{Y^0} \les$$

$$i 2^{(\frac{n}{2}-s)i} ||u_i||_{Z^s} ||v_{j, \geq 2^{j-i}}||_{Z^{s}}$$

Together with the previous estimate this entitles us to the claim

\beq \label{y12}
||\tilde{B}(u_i,v_{j})||_{X^{s,\q,1}_{j}} \les i^2 2^{(\frac{n}{2}-s)i} ||u_i||_{Z^s} ||v_{j}||_{Z^{s}}
\eeq

We obtain the same estimate for the case $j \leq i$ and $|i-k| \leq 2$. 

For the case $|i-j| \leq 2$ we can just import and then modify (\ref{y11}) to

\beq \label{y13}
||\tilde{B}(u,v)||_{X^{s,\q,1}_{k}} \leq 2^{(\frac{n}{2}+s)k} 2^{-2js} ||u||_{X^{s,\q,1}_{i}} ||v||_{X^{s,\q,1}_{j}} \les 
\eeq

$$j^2 2^{(\frac{n}{2}+s)k} 2^{-2js} ||u||_{Z^{s}_{i}} ||v||_{Z^{s}_{j}}$$

For $s > \frac{n}{2}$ we can add up the estimates (\ref{y12}) and (\ref{y13}) with respect to $i,j,k$ in (\ref{y14}) to claim

$$||\tilde{B}(u,v)||_{X^{s,\q,1}} \les C_{s} ||u||_{Z^{s}} ||v||_{Z^{s}}$$

Since $X^{s,\q,1}$ controls the $Z^s$ norm, this ends our argument. Hence we succeeded to show

$$Z^{s} \cdot Z^{s} \ra Z^{s}$$

The rest of the argument needed to complete the claim (\ref{A}) is similar to the one we provided for $X^{s,\q,1} + \bar{X}^{s,\q,1}$ with the obvious adjustments. 

\end{proof}

\begin{proof}[Proof of (\ref{M})] We derive this estimate mainly by duality and conjugation from the estimates we proved for the algebra properties. 

{\mathversion{bold} $i \leq j$ \ and \ $|k-j| \leq 2$}; \rm we can reduce the matters to the case $k=j$ and then use the equivalence:

$$(Z_i^0+\bar{Z}^0_i) \cdot W^{0}_j \ra W^0_j \Leq (Z_i^0+\bar{Z}^0_i) \cdot Z^{0}_j \ra Z^0_j$$

Therefore we obtain something similar to (\ref{y12}):

\beq \label{y39}
||\tilde{B}(u,v)||_{W^{s}_{j}} \les i^2 2^{(\frac{n}{2}-s)i} ||u||_{Z^s_i+\bar{Z}^s_i} ||v||_{W^{s}_j}
\eeq

\vspace{.1in}

{\mathversion{bold} $j \leq i$ \ and \ $|k-i| \leq 2$}; we can reduce the matters to the case $k=i$ and then use the equivalence:

$$(Z_i^0+\bar{Z}^0_i) \cdot W^{0}_j \ra W^0_i \Leq (Z_i^0+\bar{Z}^0_i) \cdot Z^{0}_i \ra Z^0_j$$

Therefore we obtain a modified version of (\ref{y13}):

\beq \label{y15}
||\tilde{B}(u,v)||_{W^{s}_{i}} \les i^2 2^{(\frac{n}{2}-s)j} ||u||_{Z^s_i+\bar{Z}^s_i} ||v||_{W^{s}_j}
\eeq

The logarithm of the high frequency (see the $i^2$ term) is present since when deriving (\ref{y13}) we had a high-high to low interaction and this was not a problem there. But it becomes problematic here. As longs as (let's say) $10 n j \geq i$ the estimate (\ref{y15}) is fine. 

We provide a different argument for the case $10 n j \leq i$ which eliminates the $i^2$ term in (\ref{y15}). We start by pointing out that we treat in a different way the interactions with $Z^s_i$ and $\bar{Z}^s_i$. 

Assume that $u_{i} \in Z^s_i$. Making use of (\ref{y2}) we estimate

$$||\tilde{B}(u_{i,d},v_{j})||_{\mathcal{Y}^0_{i, \leq 2^{i-j}}} \les 2^{\frac{(n-1)j-i}{2}} ||v_{j}||_{X^{0,\q,1}} ||u_{i,d}||_{L^{2}} \les $$

$$ 2^{\frac{(n-1)j-i}{2}} 2^{2j} d^{-\q} ||v_{j}||_{X^{0,-\q,1}} ||u_{i,d}||_{X^{0,\q}} \les d^{-\q} ||v_{j}||_{X^{0,-\q,1}} ||u_{i}||_{Z^0}$$

In the last line we have used the fact that $10 n j \leq i$. Summing up with respect to $d$ gives us:

\beq \label{y16}
||\tilde{B}(u_{i},v_{j})||_{\mathcal{Y}^0_{i, \leq 2^{i-j}}} \les ||v_{j}||_{X^{0,-\q,1}} ||u_{i}||_{Z^0}
\eeq

Making use of (\ref{y1}) it follows:

\beq \label{y17}
||\tilde{B}(u_{i, \leq 2^{i-j}},v_{j})||_{X^{0,-\q,1}_i} \leq \sum_{d} ||\tilde{B}(u_{i, \leq 2^{i-j}},v_{j})||_{X^{0,-\q}_{i,d}} \approx 
\eeq

$$\sum_{d}  d^{-\q}  ||\tilde{B}(u_{i, \leq 2^{i-j}},v_{j})||_{L^{2}_{i,d}} \les ||\tilde{B}(u_{i, \leq 2^{i-j}},v_{j})||_{L^{2}} \les $$

$$ 2^{\frac{(n-1)j-i}{2}} ||v_{j}||_{X^{0,\q,1}} ||u_{i,\leq 2^{i-j}}||_{Y^0} \les  2^{\frac{j-i}{2}} 2^{2j} ||v_{j}||_{X^{0,-\q,1}} ||u_{i}||_{Z^{0}} \les ||v_{j}||_{X^{0,-\q,1}} ||u_{i}||_{Z^0}$$

We used again the fact that $10 n j \leq i$. 

Looking back at how we derived (\ref{y11}) we notice that at the dyadic level we obtain the estimate in $X^{0,\q}$ without any logarithm

$$||\tilde{B}(u_{i,d_1},v_{j,d_2})||_{X_{i,d_3}^{0,-\q}} \les 2^{\frac{n}{2}j} ||u_{i,d_1}||_{X^{0,\q}} ||v_{j,d_2}||_{X^{0,-\q}}$$

From this, via a standard argument, we conclude with

\beq \label{y18}
||\tilde{B}(u_{i, \geq 2^{i-j}},v_j)||_{X_{i, \geq 2^{i-j}}^{0,-\q,1}} \les j 2^{\frac{n}{2}j} ||u_{i, \geq 2^{i-j}}||_{X^{0,\q,1}} ||v_{j}||_{X^{0,-\q,1}}
\eeq

Bringing together the estimates in (\ref{y16}), (\ref{y17}) and (\ref{y18}) we obtain:

\beq \label{y30}
||\tilde{B}(u_{i},v_{j})||_{W_{i}^{s}} \les j^{2} 2^{(\frac{n}{2}-s)j} ||u_{i}||_{Z^{s}} ||v_{j}||_{W^{s}}
\eeq

Next we work with $u_{i} \in \bar{Z}^s_i$. We split $u_{i}= u_{i, \leq 2^{2i-10}} + u_{i, \geq 2^{2i-9}}$. One can easily check the inclusion $\bar{Z}_{i, \geq 2^{2i-9}}^{0} \subset Z_i^0$, therefore we obtain for free 

\beq \label{y31}
||\tilde{B}(u_{i, \geq 2^{2i-9}},v_{j})||_{W_{i}^{s}} \les j^{2} 2^{(\frac{n}{2}-s)j} ||u_{i, \geq 2^{2i-9}}||_{\bar{Z}^{s}} ||v_{j}||_{W^{s}}
\eeq

Since in all the estimates we are allowed to place a conjugate on the low frequency, the argument for $Z_{i}^0 \cdot W_{j}^0 \ra W_{i}^0$ generates an argument for $Z_{i}^0 \cdot \bar{W}_{j}^0 \ra W_{i}^0$. Then, by conjugation, we obtain the estimate for $\bar{Z}_{i}^0 \cdot W_{j}^0 \ra \bar{W}_{i}^0$.  Therefore we have the estimate:

$$||\tilde{B}(u_{i, \leq 2^{2i-10}},v_{j})||_{\bar{W}_{i}^{s}} \les j^{2} 2^{(1-s)j} ||u_{i, \leq 2^{2i-10}}||_{\bar{Z}^{s}} ||v_{j}||_{W^{s}}$$

On the other hand, $\tilde{B}(u_{i, \leq 2^{2i-10}},v_{j})$ is supported in $\bar{A}_{i, \leq 2^{2i-8}}$ and we have the inclusion $\bar{W}_{i, \leq 2i-8}^0 \subset W_{i}^0$, therefore

\beq \label{y32}
||\tilde{B}(u_{i, \leq 2^{2i-10}},v_{j})||_{W_{i}^{s}} \les j^{2} 2^{(\frac{n}{2}-s)j} ||u_{i, \leq 2^{2i-10}}||_{\bar{Z}^{s}} ||v_{j}||_{W^{s}}
\eeq

From (\ref{y31}) and (\ref{y32}) we obtain the desired estimate for $\bar{Z}^{s}_{i}$, therefore, recalling (\ref{y30}), we can conclude with:

\beq \label{y41}
||\tilde{B}(u_{i},v_{j})||_{W_{i}^{s}} \les j^{2} 2^{(\frac{n}{2}-s)j} ||u_{i}||_{Z^{s}+\bar{Z}^s} ||v_{j}||_{W^{s}}
\eeq

\vspace{.1in}

{\mathversion{bold} $|i-j| \leq 2$}; here it is relevant to consider the particular case $i=j$. We have the equivalence:

$$\bar{Z}^0_i \cdot W^{0}_i \ra W^0_k \Leq Z_i^0 \cdot Z^{0}_k \ra Z^0_i$$
 
Since the last estimate had been derived in (\ref{y12}), we have:

$$||\tilde{B}(u_{i},v_{i})||_{W_{k}^{s}} \leq k^{2} 2^{(\frac{n}{2}+s)k} 2^{-2is} ||u_{i}||_{\bar{Z}^s} ||v_{i}||_{W^{s}}$$

The estimate $Z^0_i \cdot W^{0}_i \ra W^0_k$ will be derived directly. We use (\ref{b195}) to obtain:

$$||\tilde{B}(u,v)||_{X^{0,-\q,\infty}_{k}} \les 2^{(\frac{n}{2}-2)k} ||u||_{X^{0,\q,1}_{i}} ||v||_{X^{0,\q,1}_{i}}$$

For our purpose, we can modify this estimate to

$$||\tilde{B}(u,v)||_{X^{0,-\q,\infty}_{k}} \les 2^{(\frac{n}{2}-2)k+2i} ||u||_{X^{0,\q,1}_{i}} ||v||_{X^{0,-\q,1}_{i}}$$

For general $s$ this becomes:

\beq \label{y42}
||\tilde{B}(u,v)||_{Z^{s}_{k}} \leq ||\tilde{B}(u,v)||_{X^{s,-\q,1}_{k}} \les k i^{2} 2^{(\frac{n}{2}-2+s)(k-i)} 2^{(\frac{n}{2}-s)i} ||u||_{Z^{s}_{i}} ||v||_{Z^{s}_{i}}
\eeq

A standard argument sums the estimates (\ref{y39}), (\ref{y41}) and (\ref{y42}) and gives us the claim in (\ref{M}). 

\end{proof}


\begin{thebibliography}{99}

\bibitem[Be]{be} Bejenaru, I., \emph{Quadratic Nonlinear Derivative Schr\"odinger Equations - Part 1}, http://www.arxiv.org/abs/math.AP/0512041

\bibitem[Ch]{ch1} Chihara, H., \emph{Gain of regularity for semilinear Schrödinger equations}, Math. Ann 315 (1999), no. 4, 529-567

\bibitem[KaKo]{he} Koch, H. and Kato, J., \emph{Uniqueness of the modified Schr\"odinger map in $H^{\frac{3}{4}+\e}$}, http://arxiv.org/abs/math.AP/0508423

\bibitem[Ka]{ka} Kato, J., \emph{Existence and uniqueness of the solution to the modified Schr\"odinger map}, Math. Res. Lett, to appear, 2005

\bibitem[KeNa]{k1} Kenig, C.E and Nahmod, A., Personal communication, 2005

\bibitem[KePoVe2]{k2} Kenig, C.E., Ponce, G. and Vega, L, \emph{Smoothing effects and local existence theory for the generalized nonlinear {S}chr\"odinger equations}, Invent. Math., 134 (1998), no. 3, 489-545


\bibitem[Kl]{kl} Klainerman, S., \emph{The null condition and global existence to nonlinear to nonlinear wave equations}, Nonlinear systems of PDE in applied mathematics, Part 1 (Santa Fe, N.M., 1984), Lectures in Appl. Math., vol 23, Amer. Math. Soc., Providence, RI, 1986, MR 87h:35217   



\bibitem[NaStUh]{n1} Nahmod, A., Stefanov, A., Uhlenbeck, K., \emph{On Schroedinger maps}, CPAM 56 (2003), 114-151

\bibitem[Tao]{tao1} Tao, T., \emph{Global regularity of wave maps II. 
 Small energy in two dimensions}, Comm. Math. Phys. 224 (2001), 443-544

\bibitem[Ta1]{ta1} Tataru, D., \emph{Rough Solutions for the Wave Maps Equation}, Amer. J. Math., 127(2005), no. 2, 293-377 

\bibitem[Ta2]{ta2} Tataru, D., \emph{The Wave Maps Equation}, Bull. of the AMS, vol. 41(2004), no. 2, 185-204


\end{thebibliography}
\end{document}